\documentclass[11pt,a4paper,leqno]{article}
\usepackage{amsmath}
\usepackage{amssymb}
\usepackage{graphicx}

\advance\topmargin by -2.2cm
\newcommand{\barr}{\begin{array}}
\newcommand{\earr}{\end{array}}
\newcommand{\beqq}{\begin{equation}}
\newcommand{\eeqq}{\end{equation}}
\newcommand{\beao}{\begin{eqnarray*}}
\newcommand{\eeao}{\end{eqnarray*}\noindent}
\newcommand{\beam}{\begin{eqnarray}}
\newcommand{\eeam}{\end{eqnarray}\noindent}

\newcommand{\si}{\sigma}
\newcommand{\al}{\alpha}

\newcommand{\vep}{\varepsilon}

\newcommand{\wh}{\widehat}

\newcommand{\lra}{\longrightarrow}

\begin{document}

{\large\bf To which extent is the membrane potential in a neuron between successive spikes adequately modelled 
by a (continuous) semimartingale?}\\   

Reinhard H\"opfner, Universit\"at Mainz \\ 

{\small 
{\bf Abstract: }
We consider $p$-variations in some membrane potential data --viewed as a function of the step size in case where $p$ is fixed, or viewed as a function of $p$ in case where the step size is fixed-- and compare their shape with results in Jacod and Ait-Sahalia \cite{AJ-09a} which do hold for general semimartingales. We obtain the following conclusion: 
in non- or very rarely-spiking cases the membrane potential behaves as a semimartingale, in some cases as a semimartingale with jumps. Once the neuron is spiking, a semimartingale modelization is no longer adequate for the membrane potential between successive spikes, even if interspike intervals are relatively long. 


}

\vskip1.5cm


We take a new look on two data sets recording the membrane potential in a pyramidal neuron (intracelluar recording) which belongs to a cortical slice observed in vitro (representing an active network) under different experimental conditions. The neuron under observation receives synaptic input from a large number of other neurons in the slice. Stimulating the slice --and thus the networking properties of all neurons belonging to the slice--  by a potassium bath, W.\ Kilb (Institute of Physiology, University of Mainz) recorded 'Zelle~3' in 2004 and '17Sept08{\_}023' in 2008. The data are shown in figures \ref{fig19} and \ref{fig20}. 
In 'Zelle~3', 10 different concentrations of potassium correspond to 10 different data sets (called 'levels' below, obtained under 3, 4, 5, 6, 7, 8, 9, 10, 12, 15 mM of K) observed over 60 seconds each. In '17Sept08{\_}023', one potassium level (5 mM of K) was kept constant over a much longer time interval; in the present note, for ease of comparison with 'Zelle~3', we use only the first 60 seconds of observation from this data set.

\begin{figure}[t]
\vskip0.8cm
\begin{center}
\includegraphics[width=12.5cm]{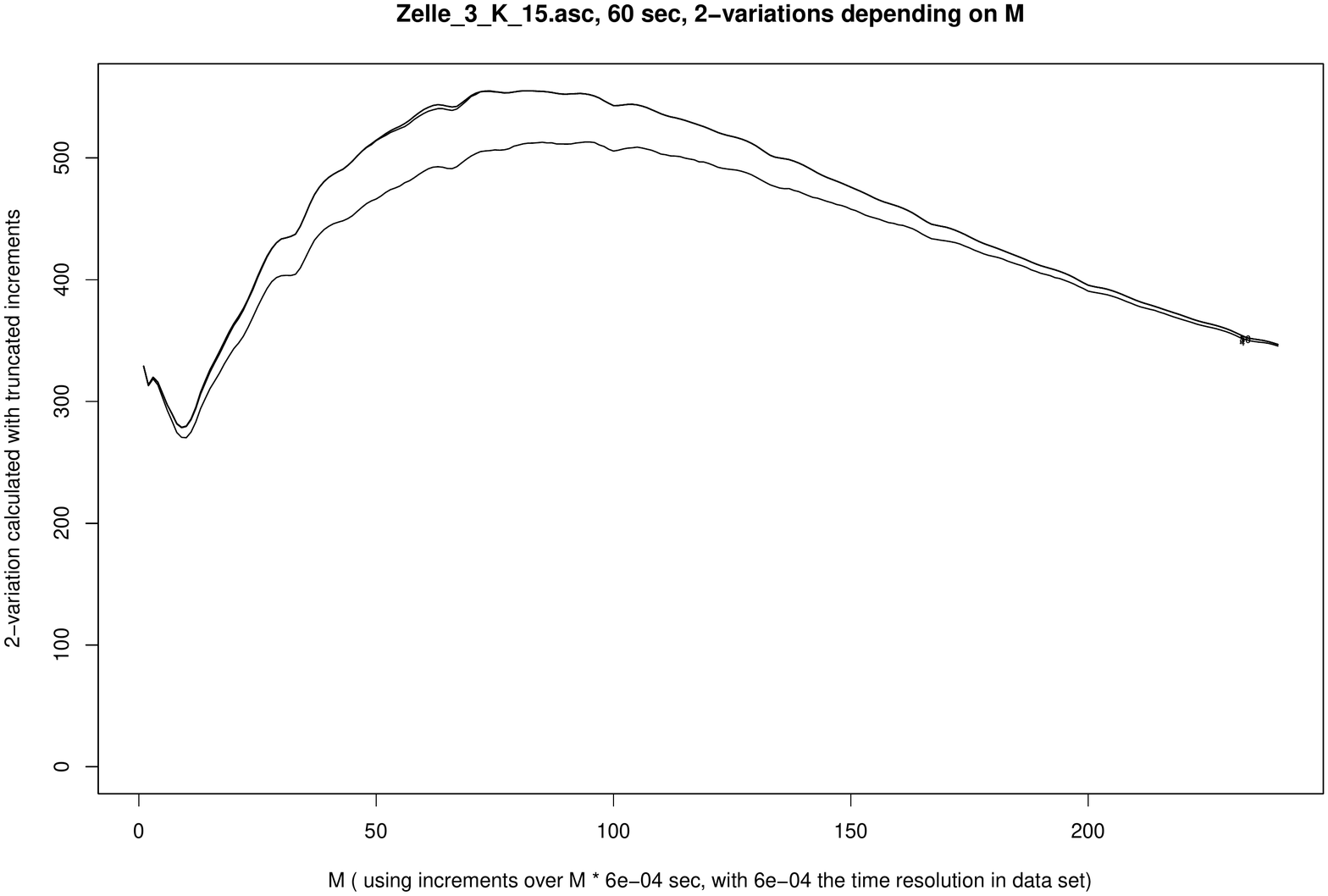} 
\caption{\small Membrane potential 'Zelle~3' level~10 (15 mM of K): plotting truncated 2-variations $M\to V_\Gamma(2,\Delta,M)$ as  defined by (\ref{mypvariationsegment}) for $1\le M\le 240$, in increasing order for $\Gamma = 4,8,10,16$; no further changes above $\Gamma=10$.}
\label{fig9}
\end{center}
\end{figure}

In terms of a diffusion process modelization, 'Zelle~3' has been considered in \cite{H-07}, '17Sept08{\_}023' in \cite{Jahn-09} (section 5.3 there, using the estimation method of \cite{H-07}). 
In \cite{H-07} and \cite{Jahn-09}, assuming that the membrane potential between successive spikes (more precisely: sufficiently away from the spikes) can be modelled as a time homogeneous diffusion process, nonparametric estimates for diffusion coefficient and drift made appear a linear mean-reverting drift combined with either a constant or a linear or a 'bowl-shaped' diffusion coefficient: these cases correspond to Ornstein-Uhlenbeck~(OU) type, Cox-Ingersoll-Ross~(CIR) type or --in the language of \cite{FS-08}-- Pearson~(P) type diffusions. Analyzing the data in the same way in smaller time windows, we can assert that the assumption of time homogeneity seems  well satisfied in '17Sept08{\_}023' (here W.\ Kilb had used a new type of electronic stabilization device), and reasonably well satisfied in several potassium levels of 'Zelle~3'; obvious exceptions are the 'low' levels 1, 2, 4 where strong time inhomogeneities appear (discussed for level~1 in \cite{H-07}, section 4.6).

In this note, we consider for $p{=}2$ or $p{=}4$ fixed $p$-variations in the membrane potential data as a function of the step size, i.e.\ the length of the time intervals over which variations are calculated, and compare these to simulated diffusion equivalents whose drift and diffusion coefficients are as estimated in \cite{H-07} or \cite{Jahn-09}. 
Relying on recent results of Ait-Sahalia and Jacod \cite{AJ-09a}, we then ask the question to which extent a (continuous) semimartingale model is in fact adequate for the membrane potential between successive spikes. 
In our data, a surprising difference appears --and in particular in the same neuron 'Zelle~3'-- between spiking and non-spiking regimes. 
We then fix the step size and consider $p$-variations in our data as a function of $p$: again the same striking difference between spiking and non-spiking regimes arises. When spiking is sufficiently frequent 
\begin{figure}[p]
\begin{center}
\includegraphics[width=12.5cm]{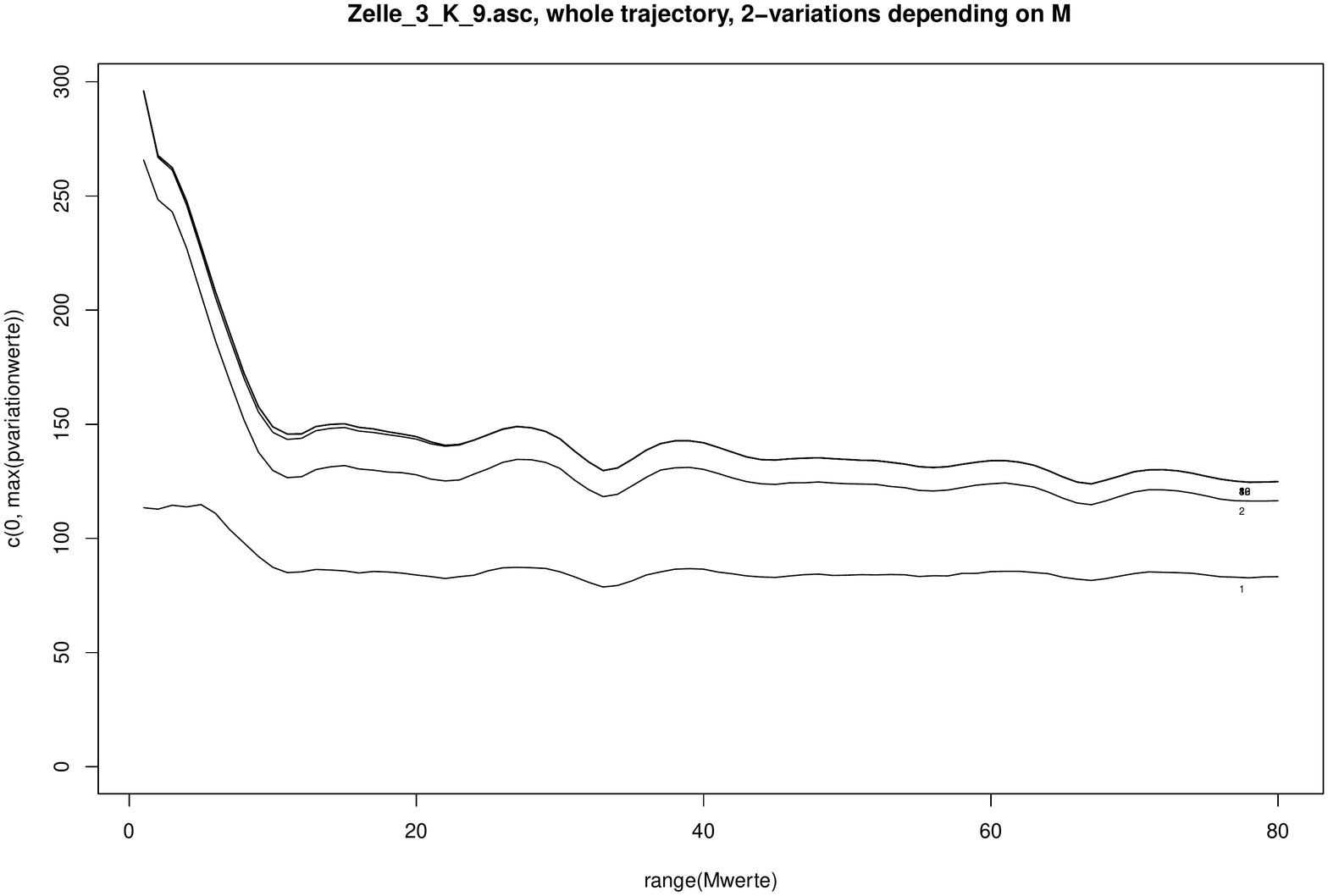}
\caption{\small Membrane potential Zelle~3 level~7 (9 mM of K): plotting truncated 2-variations $M\to V_\Gamma(2,\Delta,M)$ as defined by (\ref{mypvariationsegment}), in increasing order for $\Gamma \in\{1,2,4,8,10,16,32\}$; no further changes above $\Gamma= 8$.}
\label{fig1}
\includegraphics[width=12.5cm]{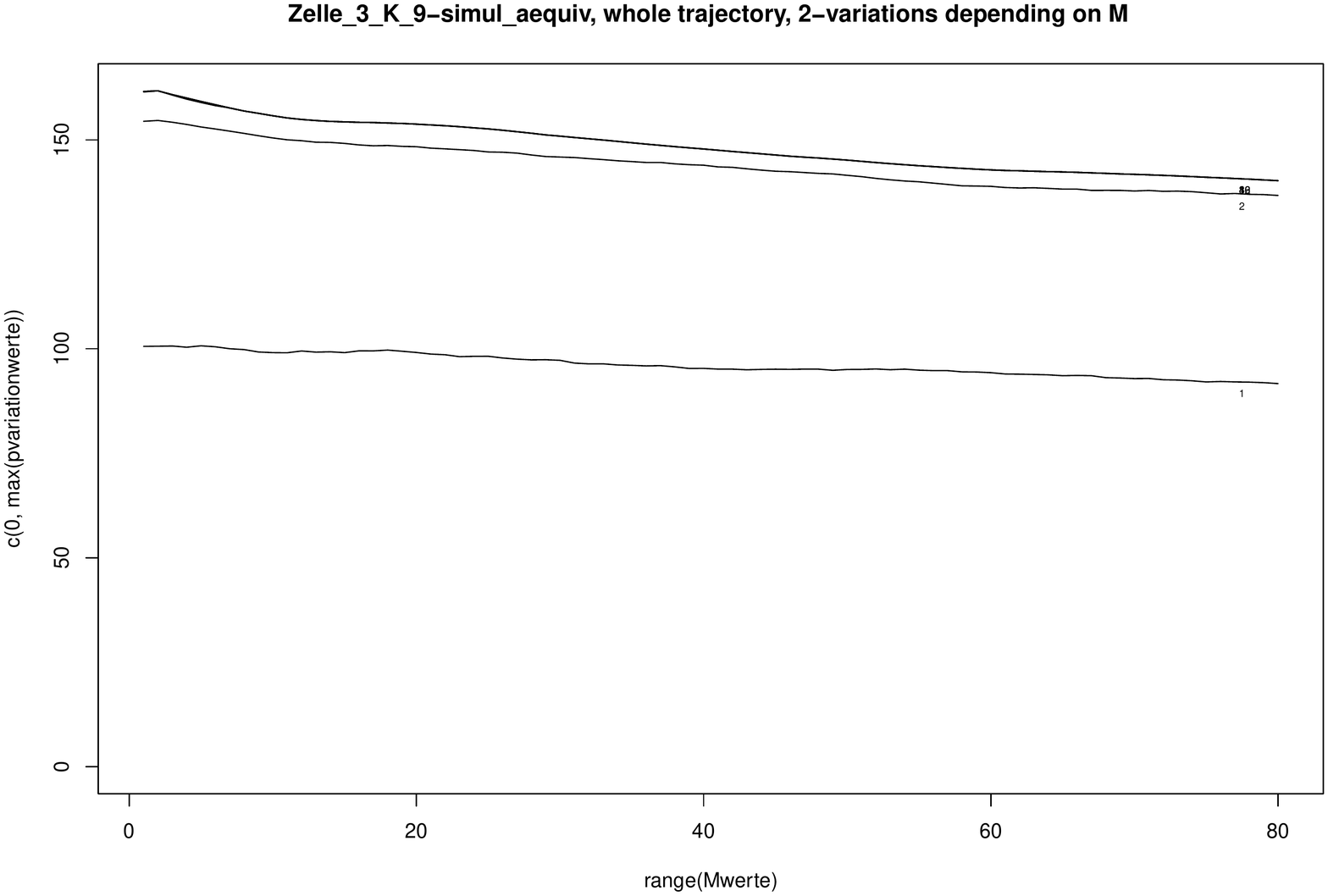}
\caption{\small Simulated P type diffusion equivalent (\cite{H-07}, section 4.4) for Zelle~3 level~7: 2-variations plotted in analogy to figure \ref{fig1}.}
\label{fig2}
\end{center}
\end{figure}
\begin{figure}[p]
\begin{center}
\includegraphics[width=12.5cm]{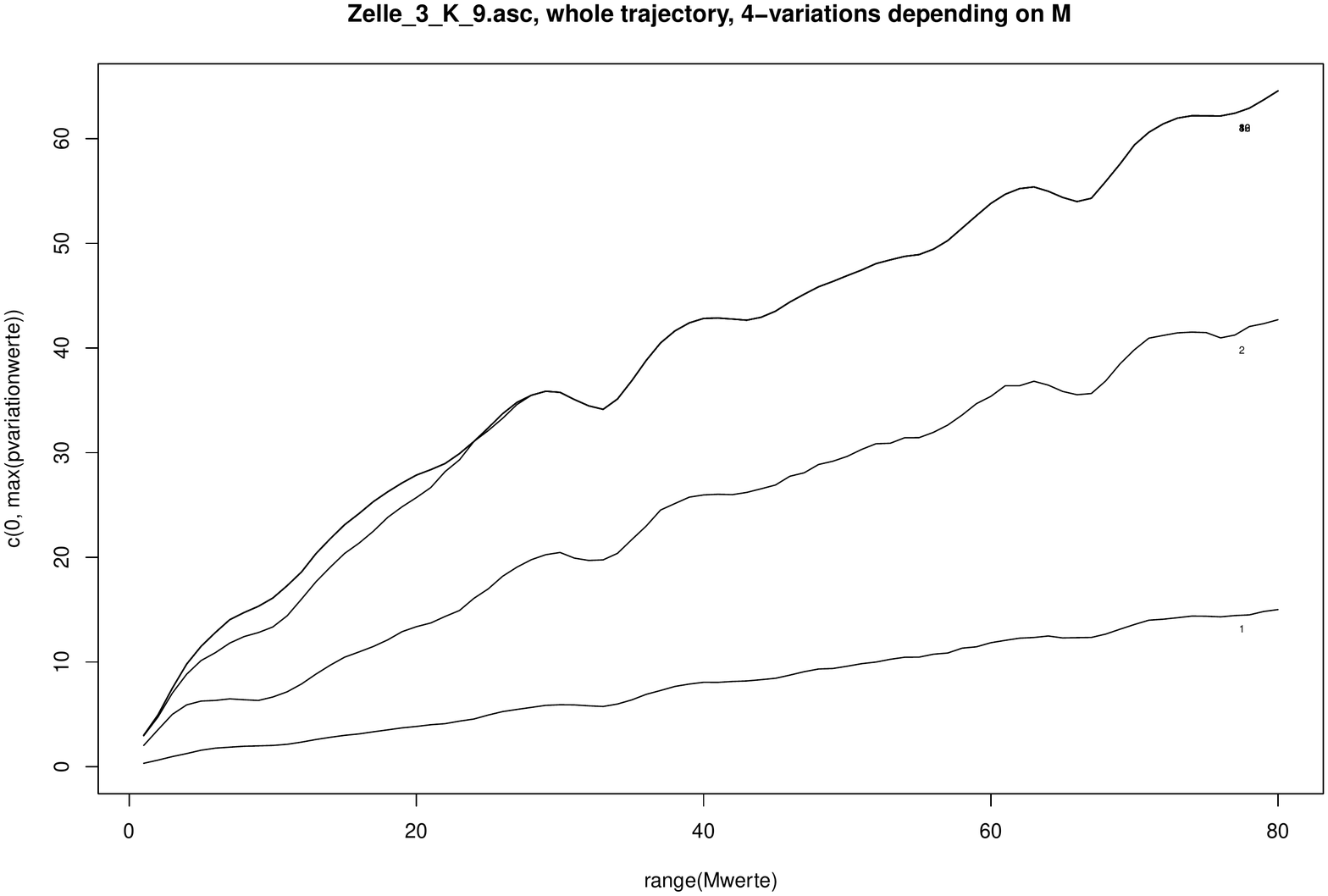}
\caption{\small Membrane potential Zelle~3 level~7 (9 mM of K): plotting truncated 4-variations $M\to V_\Gamma(4,\Delta,M)$ as defined by  (\ref{mypvariationsegment}), in increasing order for $\Gamma \in\{1,2,4,8,10,16,32\}$; no further changes above  $\Gamma= 8$.}
\label{fig3}
\includegraphics[width=12.5cm]{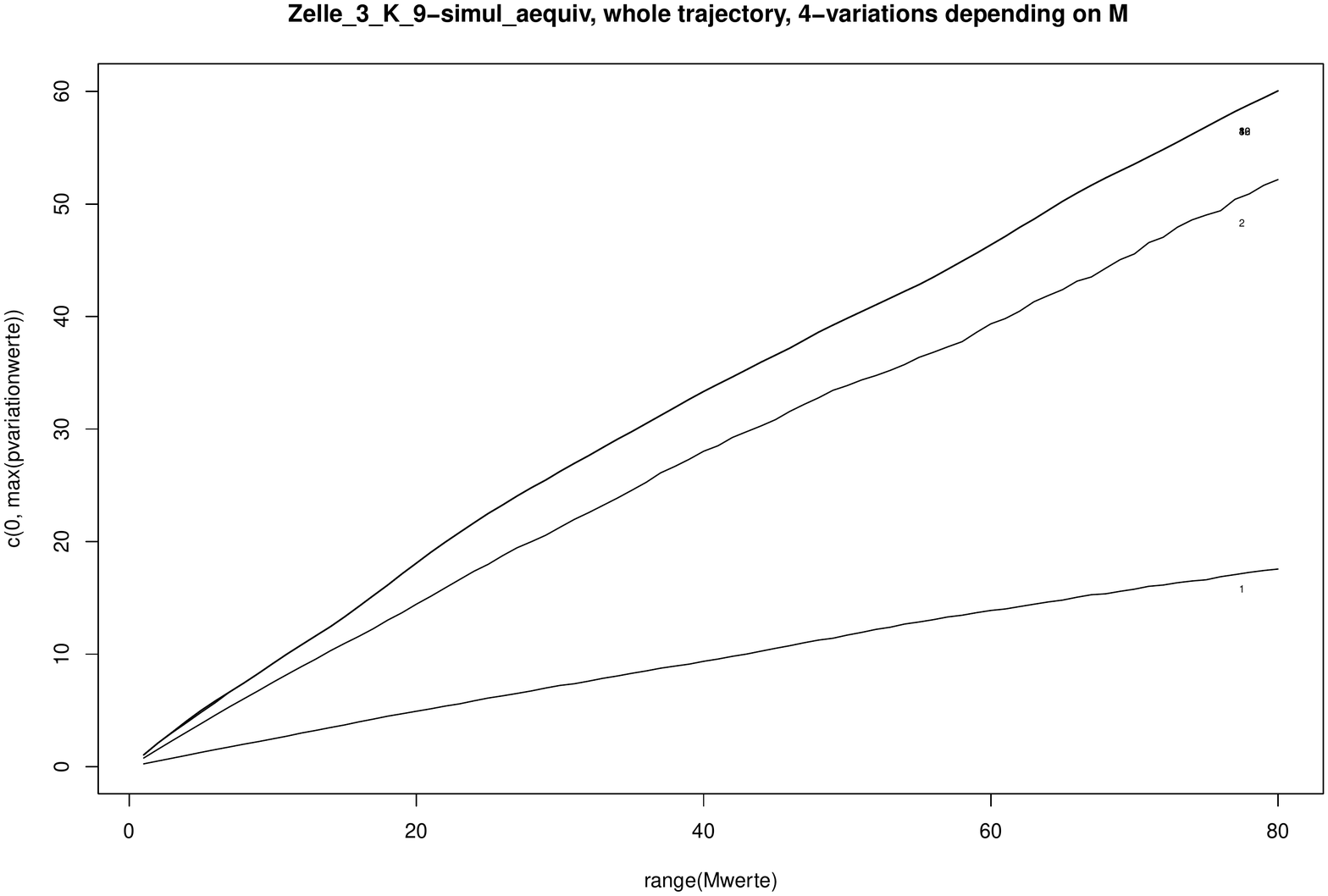} 
\caption{\small Simulated P type diffusion equivalent for 'Zelle~3' level~7 (\cite{H-07}, section 4.4): 4-variations plotted in analogy to figure \ref{fig3}.}
\label{fig4}
\end{center}
\end{figure}
(in our data, in levels 9 and 10 of 'Zelle~3', and in '17Sept08{\_}023'), pictures of power variations arise which do not agree with what we should see in a semimartingale, continuous or not, according to \cite{AJ-09a}. This is not simply an effect of noisy observation (equally observable in all levels of 'Zelle~3', but almost absent in '17Sept08{\_}023' where a different type of electrode had been used) but concerns the shape of the curve of power variations as a function of the step size when $p$ is fixed, or as a function of $p$ when the step size is fixed. In sharp contrast to this, in non- or rarely spiking regimes (the spikeless levels 1--7  of 'Zelle~3' (3--9 mM of K), and level 8 (10 mM of K) where one single spike is emitted during the overall observation time of 60 seconds) the pictures of power variations agree very well with what is to be expected for a semimartingale --up to secondary effects like noisy observation or feedback effects in the slice-- and with what can be seen in simulated diffusions or jump diffusions.

Our conclusion is that a semimartingale model seems adequate for neurons in non-spiking or rarely-spiking regimes, whereas something essentially different --not well captured by semimartingale modelization-- seems to prevail in spiking regimes. 

This note is organized as follows. 
Section 1 considers for fixed $p$ ($p=2$ or $p=4$) $\;p$-variations as a function of the step size over which we calculate the increments. Subsection 1.1 explains the truncated power variations which we use in this note, subsection 1.2 considers the spikeless or very rarely spiking levels of 'Zelle~3', subsection 1.3 the spiking levels of 'Zelle~3' and the frequently spiking neuron  '17Sept08{\_}023'. Section 2 considers $p$-variations as a function of $p$ for fixed step size, with an analogous program. 
I would like to stress that this note is a 'not really mathematical' paper (no theorem, no rigorous proofs, some merely 'plausible' approximations): its aim is to analyze a set of neuronal data in the light of theorems in Ait-Sahalia and Jacod \cite{AJ-09a} which do hold for very general semimartingales, and to show that some essential difference exists between spiking and non-spiking regimes (in the same neuron) in view of semimartingale modelization. I would like to thank H.\ Luhmann and W.\ Kilb for the data, and J.\ Jacod for some longer discussions on this problem.

\section*{1. Fixing the power $p$ and varying the size of the increments} 

The structure of the data is as follows. 
\begin{figure}[t]
\begin{center}
\includegraphics[width=12.5cm]{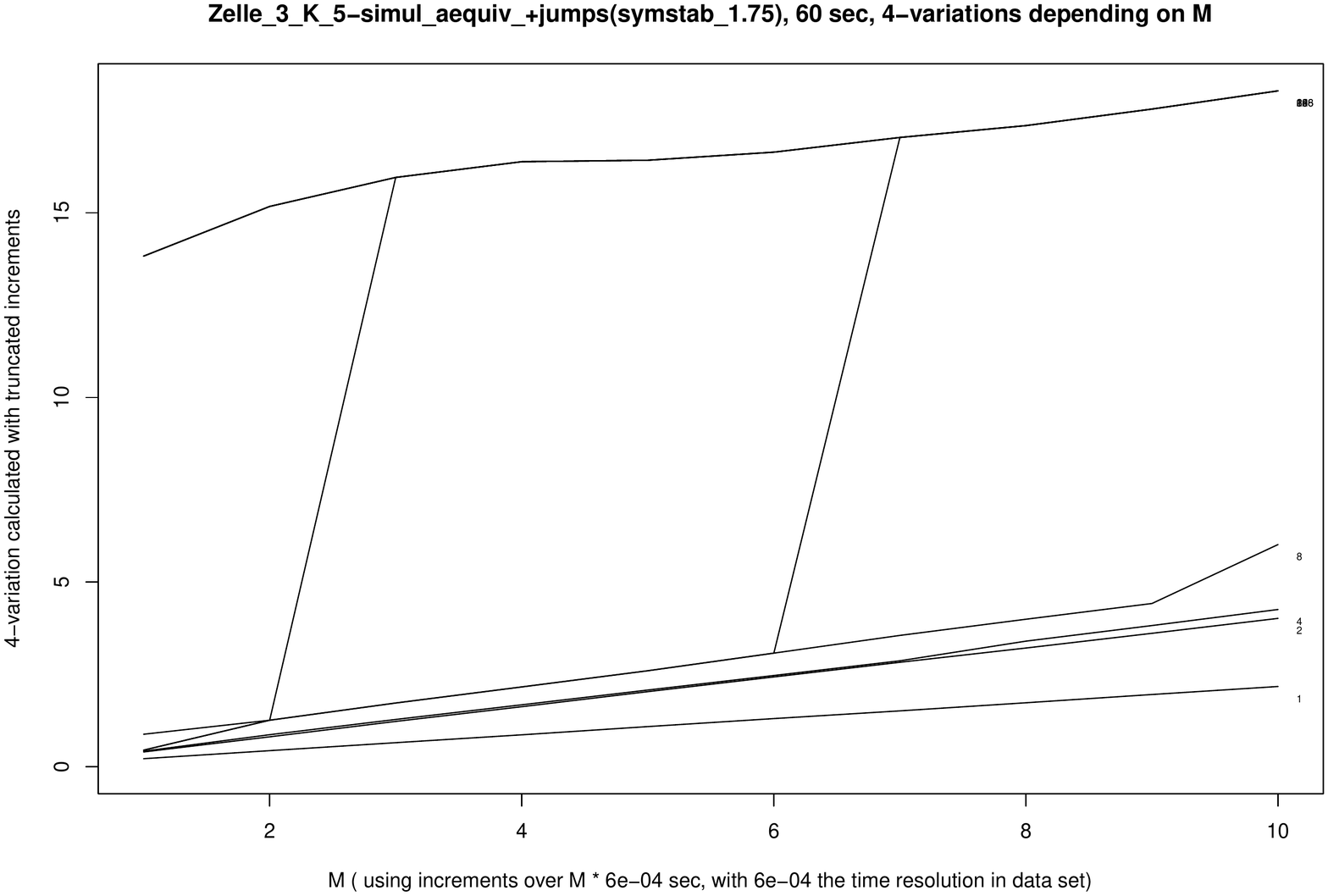}
\caption{
\small Replacing $dW_t$ in the simulated diffusion equivalent to 'Zelle~3' level~5 (\cite{H-07}, section 4.1) by $dW_t + \vep dS^\al_t$ for $\alpha=1.75$ and $\vep=0.1$: truncated 4-variations $M\to V_\Gamma(4,\Delta,M)$ as defined in (\ref{mypvariationsegment}) for $1\le M\le 10$, in increasing order for $\Gamma \in\{1,2,4,8,10,16,32, 64, 128, 256\}$; no further changes above  $\Gamma=64$.}
\label{figeffectsofjumps}
\end{center}
\end{figure}
The different experiments in 'Zelle~3' record membrane potentials at times $t_i:= i\Delta$, $0\le i\le 100001$, with $\Delta=6\cdot 10^{-4}$ [sec], thus with total observation time $T=60$ [sec]. The data  '17Sept08{\_}023' are considered only in restriction to the first 60 seconds of observation: here the time grid is $t_i=i\Delta$ with $\Delta=2 \cdot 10^{-4}$ [sec], $0\le i\le 300001$. The  measurement is in millivolt [mV], formally with three decimals, but with an information that the third decimal is not reliable at all. When spikes are present in the data set, we remove time neighbourhoods $(\tau-0.12,\tau+0.18)$ [sec] centred at the spike times $\tau$ from the data, in order to exclude any influence of the typical shape of the spike. We  calculate increments, $p$-variations, $\ldots$ over spikeless time segments, and then add up corresponding terms coming from different segments. The estimators which we use for drift and diffusion coefficient in a discretely observed diffusion are those of \cite{H-07}, and are not explained here. Here we explain the way we calculate a $p$-variation over a time segment, and over a collection of time segments, in view of application to the membrane potential in a neuron which can emit spikes. In the present section, we concentrate on fixed power $p$ and vary the step size, i.e.\ the length of the time intervals on which increments are evaluated.

\subsection*{1.1. Truncated $p$-variations for neuronal data} 
First, for varying choices of a truncation factor $0<\Gamma <\infty$ and for multiples $M$ of the step size $\Delta$ prescribed by the data, for $p\ge 2$ fixed, we define  
\beqq\label{mypvariation}
V_{t_0,t_1, \Gamma} (p,\Delta,M) \quad:=\quad 
\frac1M\; \sum_{i=i_0}^{i_1-M} \left| X_{(i+M)\Delta} - X_{i\Delta}  \right|^p \; 
1_{ \{ \;\left| X_{(i+M)\Delta} - X_{i\Delta}  \right| \;\le\; 3\, \sqrt{\Delta M}\; \Gamma\; \} }  
\eeqq
with respect to one spikeless segment $[t_0,t_1] = [i_0\Delta,i_1\Delta]$  (we define a spikeless segment as a maximal interval between $\tau_{r-1}+0.18$ and $\tau_r-0.12$, avoiding neighbourhoods of the successive spike times $\tau_{r-1}$, $\tau_r$  as defined above). Second, based on (\ref{mypvariation}), we define 
\beqq\label{mypvariationsegment}
V_{\Gamma} (p,\Delta,M) \quad:=\quad   
V_{t_{0,1},t_{1,1}, \Gamma} (p,\Delta,M)  + \ldots + V_{t_{0,\ell},t_{1,\ell}, \Gamma} (p,\Delta,M)
\eeqq 
for the whole membrane potential trajectory up to time $T=60$ [sec], where $[t_{0,1},t_{1,1}], \ldots, [t_{0,\ell},t_{1,\ell}]$ denotes the collection of spikeless segments (including an initial $[0,\tau_1-0.12]$ before the first spike and a final $[\tau_{\ell-1}+0.18,T]$ after the last spike). For spikeless membrane potentials or for the simulated diffusion equivalents, the full interval $[0,T]$ is the unique segment.

With truncation factor $\Gamma$ increasing to $\infty$ in (\ref{mypvariationsegment}), we will finally capture all jumps of a semimartingale trajectory up to time $T$, or all increments from spikeless segments in a membrane potential data set, and will arrive for $\Gamma$ tending to $\infty$ at 
\beqq\label{eq-1}
V_{t_0,t_1} (p,\Delta,M) \;:=\; V_{t_0,t_1, \infty} (p,\Delta,M)  =
\frac1M\; \sum_{i=i_0}^{i_1-M} \left| X_{(i+M)\Delta} - X_{i\Delta}  \right|^p  
\eeqq 
for a single spikeless segment $[t_0,t_1]$, and at 
\beqq\label{eq-2}
V (p,\Delta,M) \;:=\;  
V_{t_{0,1},t_{1,1}} (p,\Delta,M)  + \ldots + V_{t_{0,\ell},t_{1,\ell}} (p,\Delta,M) \;. 
\eeqq
for the whole membrane potential trajectory. 
For $t_0\le s_0 < t_1$, consider $M\Delta$-step $p$-variations on $[s_0,t_1]$ as defined in Ait-Sahalia and Jacod (\cite{AJ-09a}, formula (9)): 
$$
\wh B_{s_0,t_1}(p,\Delta,M) \;:=\;  \sum_{k=1}^{ \lfloor \frac{t_1-s_0}{M\Delta} \rfloor }
\left| X_{s_0+kM\Delta} - X_{s_0+(k-1)M\Delta}  \right|^p   
$$
and note that the right hand side of (\ref{eq-1}) equals 
\beqq\label{eq-3}
\frac1M\; \sum_{j=0}^{M-1} \wh B_{(i_0+j)\Delta\,,\,t_1}(p,\Delta,M) \;. 
\eeqq   
Averaging over $j=0,1,\ldots, M{-}1$ in (\ref{eq-3}) allows to make use of all $M\Delta$--step increments available in the time window $[t_0,t_1]$. {\em Heuristically}, for $M\Delta$ sufficiently small, all summands in (\ref{eq-3}) should be very close to $\wh B_{t_0,t_1}(p,\Delta,M)$, thus we {\em will make the following approximation during the present note} :    
\beqq\label{eq-4}
\frac1M\; \sum_{j=0}^{M-1} \wh B_{(i_0+j)\Delta\,,\,t_1}(p,\Delta,M)  \;\;\approx\;\; \wh B_{t_0,t_1}(p,\Delta,M)  \;. 
\eeqq

\begin{figure}[p]
\begin{center}
\includegraphics[width=12.5cm]{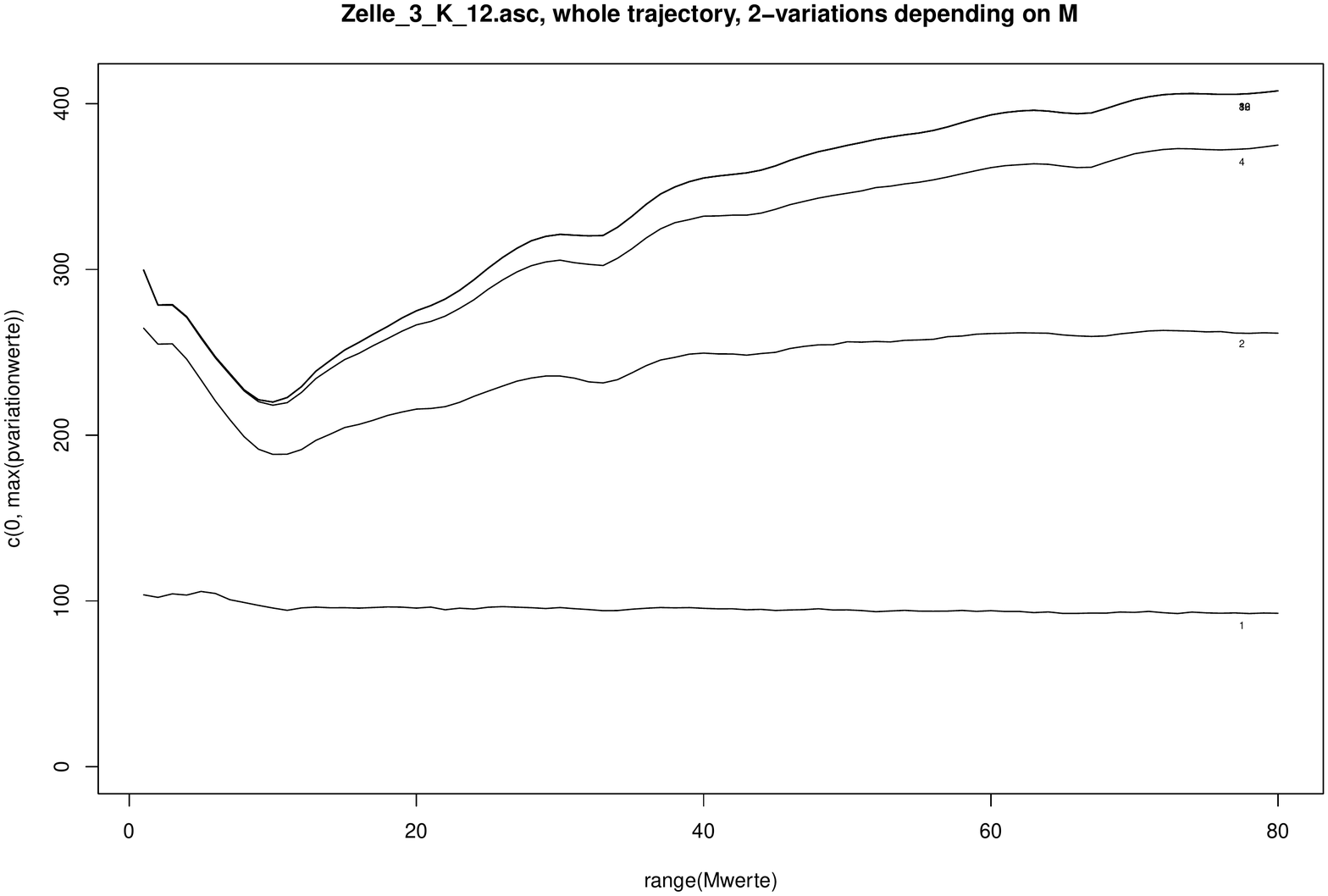}
\caption{\small Membrane potential 'Zelle~3' level~9 (12 mM of K): plotting the 2-variations $M\to V_\Gamma(4,\Delta,M)$ as defined by (\ref{mypvariationsegment}), in increasing order for $\Gamma \in\{1,2,4,8,10,16,32\}$; no further changes above  $\Gamma=8$.}
\label{fig5}
\includegraphics[width=12.5cm]{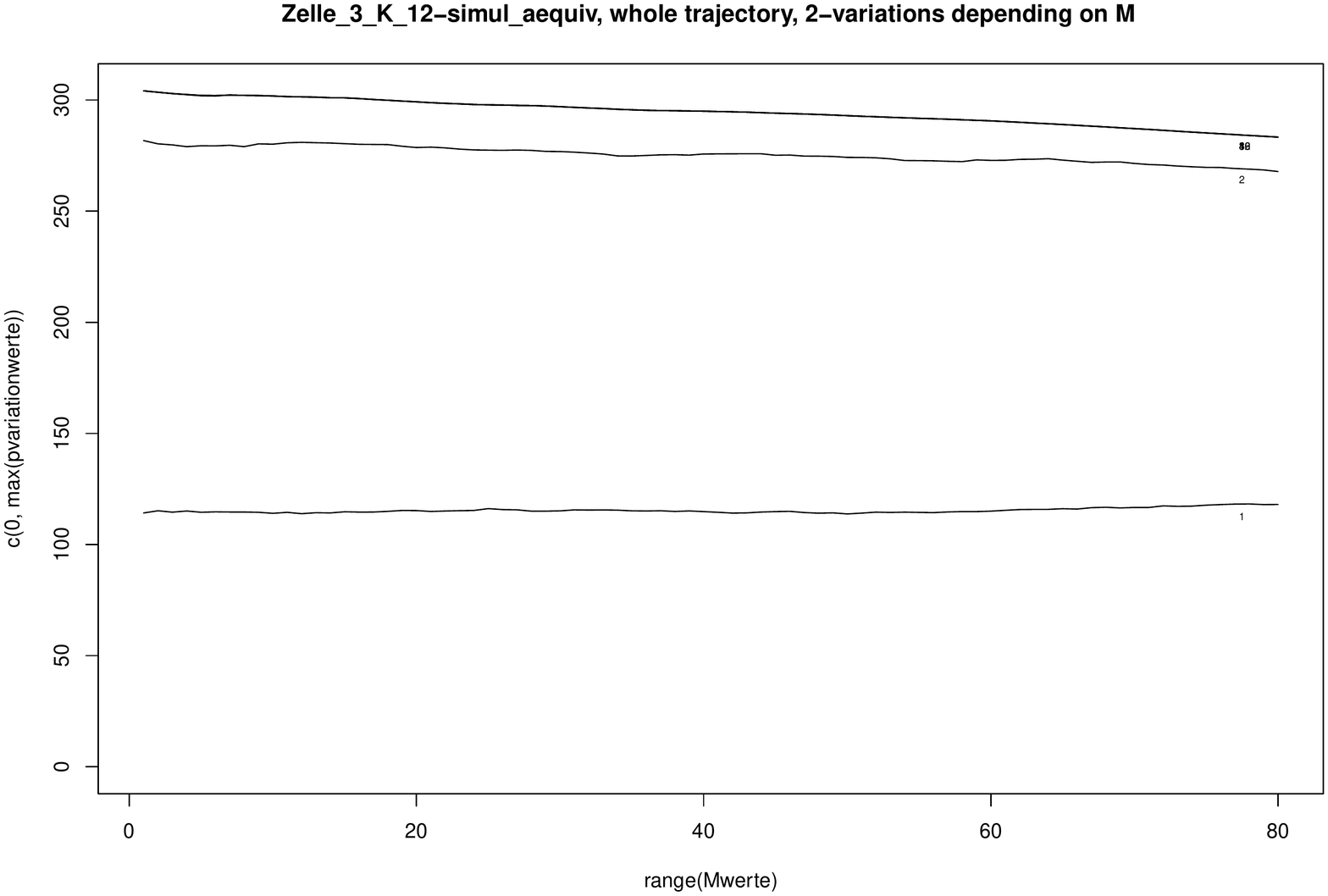} 
\caption{\small Simulated CIR type diffusion equivalent  for 'Zelle~3' level~9 (\cite{H-07}, section 3.3): 2-variations plotted in analogy to figure \ref{fig5}.}
\label{fig6}
\end{center}
\end{figure}
\begin{figure}[p]
\begin{center}
\includegraphics[width=12.5cm]{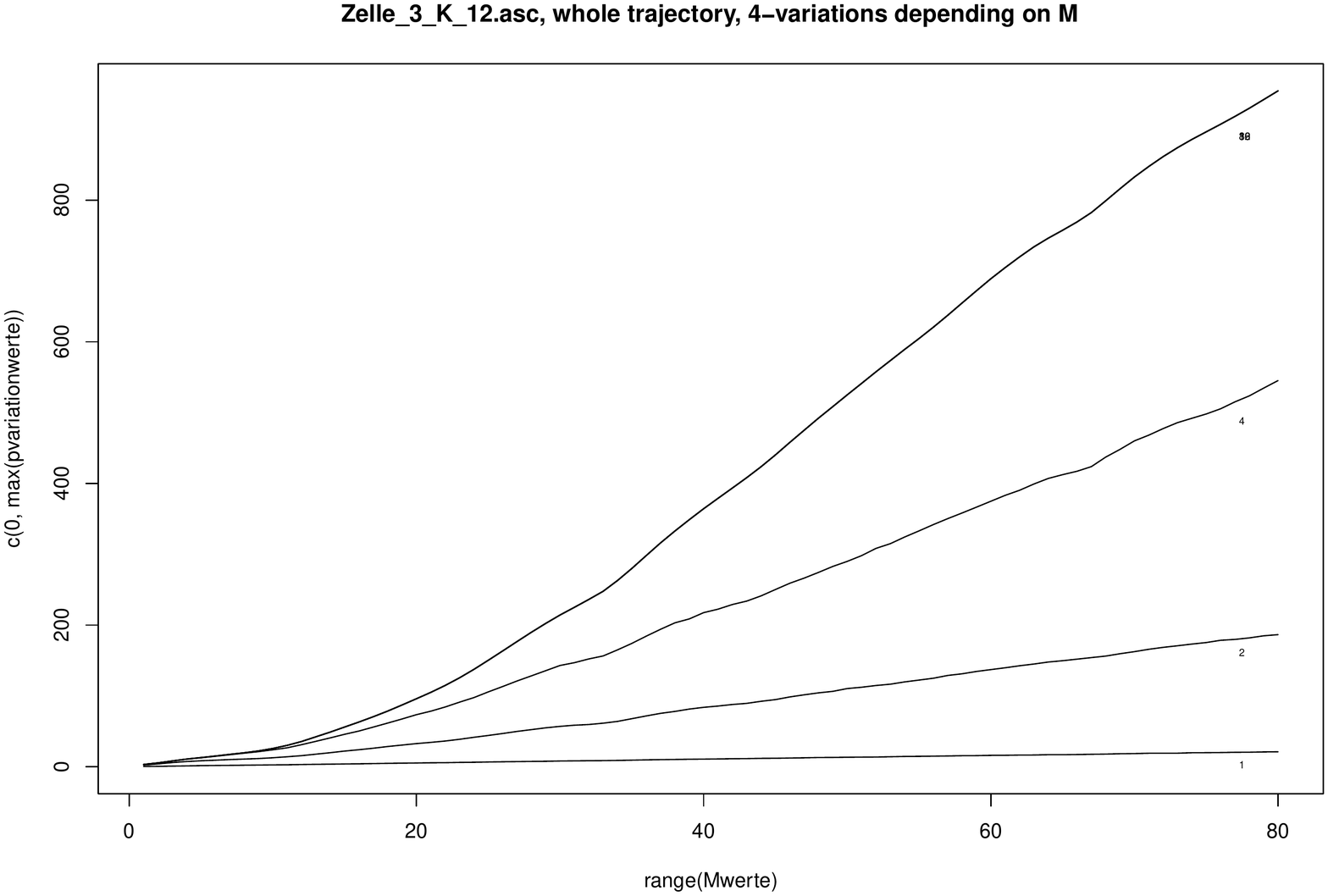}
\caption{\small Membrane potential 'Zelle~3' level~9 (12 mM of K): plotting the 4-variations $M\to V_\Gamma(4,\Delta,M)$ as defined by (\ref{mypvariationsegment}), in increasing order for $\Gamma \in\{1,2,4,8,10,16,32\}$; no further changes above  $\Gamma=8$. For $\Gamma\ge 8$, we calculate values $8.17$ for $M=3$, $5.42$ for $M=2$, $3.20$ for $M=1$. }
\label{fig7}
\includegraphics[width=12.5cm]{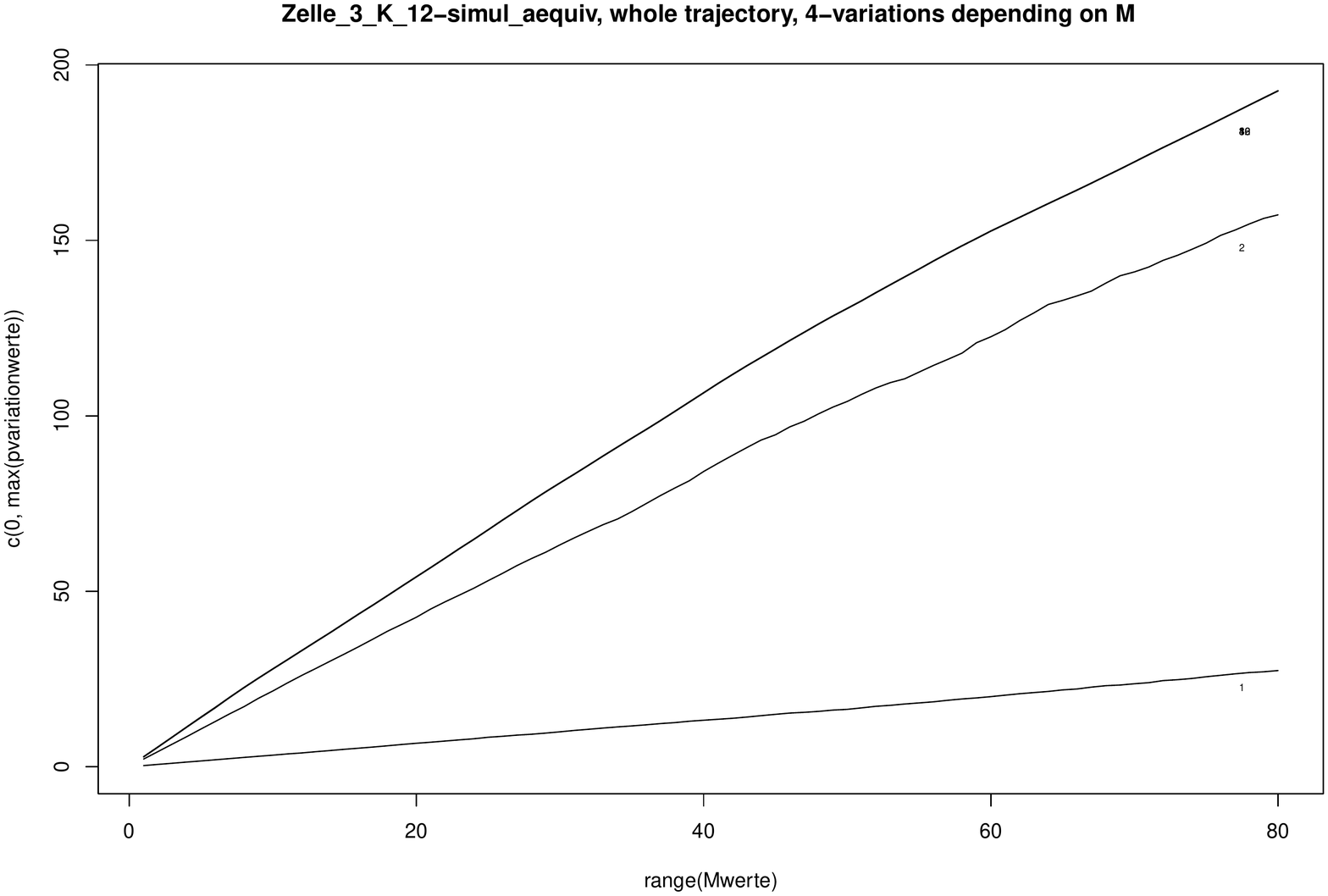}
\caption{\small Simulated CIR type diffusion equivalent for 'Zelle~3' level~9 (\cite{H-07}, section 3.3): 4-variations plotted in analogy to figure \ref{fig7}.}
\label{fig8}
\end{center}
\end{figure}

Since the time resolution $\Delta$ of the data set cannot be modified by the statistician, asymptotic results as given in Jacod and Ait-Sahalia (\cite{AJ-09a}, (11)--(13)) for time-step tending to $0$ have to be mimicked through variation of multiples $M$ of $\Delta$. Assuming that $M\Delta$ is sufficiently small for the $M$ which we consider, and that the data in fact do correspond to a discretely observed semimartingale $\xi=(\xi_s)_{s\ge 0}$, we may read  (\cite{AJ-09a}, (11)+(7)+(10)) as
\beam\label{eq-5}
\wh B_{t_0,t_1}(4,\Delta,M) &\approx& B_{t_0,t_1}(4) \,:=  \sum_{t_0\le s\le t_1}|\Delta\xi_s|^4
\quad\mbox{as $M = \ldots, 3,2,1$ gets small}
\\ \label{eq-6}
\wh B_{t_0,t_1}(2,\Delta,M) &\approx& A_{t_0,t_1}(2) + B_{t_0,t_1}(2) \,:= \int_{t_0}^{t_1}\si^2(\xi_s)\,ds +  \sum_{t_0\le s\le t_1}|\Delta\xi_s|^2 
\quad\mbox{as $M$ gets small}
\\ \label{eq-7}
\mbox{for $\xi$ continuous}&:& 
\wh B_{t_0,t_1}(4,\Delta,M) \;\approx\;  M\cdot  \wh B_{t_0,t_1}(4,\Delta,1) 
\quad\mbox{as $M$ gets small}   \;. 
\eeam
By (\ref{eq-5}),  $\,4$-variations stabilizing at a strictly positive 'limit' when $M$ gets small indicate presence of jumps in the semimartingale $\xi$. For $\xi$ continuous, $\,4$-variations should be linear in $M$ as long as $M$ is small, as a consequence of (\ref{eq-7}). On every segment, by (\cite{AJ-09a}, theorem 1), this is a dichotomy which represents a test for presence of jumps in a semimartingale $\xi=(\xi_t)_{t\ge 0}$ recorded at time resolution $\Delta$. Putting together the segments as in (\ref{eq-2}) above, we rephrase the test of \cite{AJ-09a} in the following form:  as $M=\ldots 3,2,1$ gets small, 
\beam\label{eq-5neu}
\mbox{for $\xi$ with jumps}&:& V(4,\Delta,M) \;\;\mbox{stabilizes at a strictly positive 'limit'} \;; 
\\ \label{eq-7neu}
\mbox{for $\xi$ continuous}&:&   M \;\lra\; V(4,\Delta,M) \quad\mbox{is linear}  \;. 
\eeam

The results of Ait-Sahalia and Jacod \cite{AJ-09a} being asymptotic results for shrinking time grids on which the process is observed,  reformulations such as (\ref{eq-5})+(\ref{eq-7}) or (\ref{eq-5neu})+(\ref{eq-7neu}) of this test hinge on the assumption that  $M\Delta$ be 'sufficiently small' for the $M$ which we wish to consider. In the data, we can not modify the time resolution $\Delta$. It may well happen that considering $M\Delta$ for $1\le M\le 5$ (say), we are not yet 'sufficiently small' in the sense of \cite{AJ-09a}. As an example, replace in the OU diffusion equivalent to 'Zelle~3' level 3 (as in \cite{H-07}, section 4.1) the driving $dW_t$ by $dW_t + \vep dS^\al_t$ for small $\vep$ where $S^\al = (S^\al_t)_{t\ge 0}$ is a symmetric stable process with index $\al\in(0,2)$. Simulating increments of $S^\al$ using Chambers, Mallows and Stuck \cite{CMS-76}, with $\Delta$ the time resolution of 'Zelle~3', the test (\ref{eq-5})+(\ref{eq-7}) will be unable to detect presence of jumps in the simulated jump diffusion for $\al$ very close to $2$, whereas in case $\al=1.75$, the jumps are detected (see figure \ref{figeffectsofjumps}) by inspection of $4$-variations for $M\le 5$.

\subsection*{1.2. Application to the non- or very rarely spiking levels of 'Zelle~3'}  
\begin{figure}[p]
\begin{center}
\includegraphics[width=12.5cm]{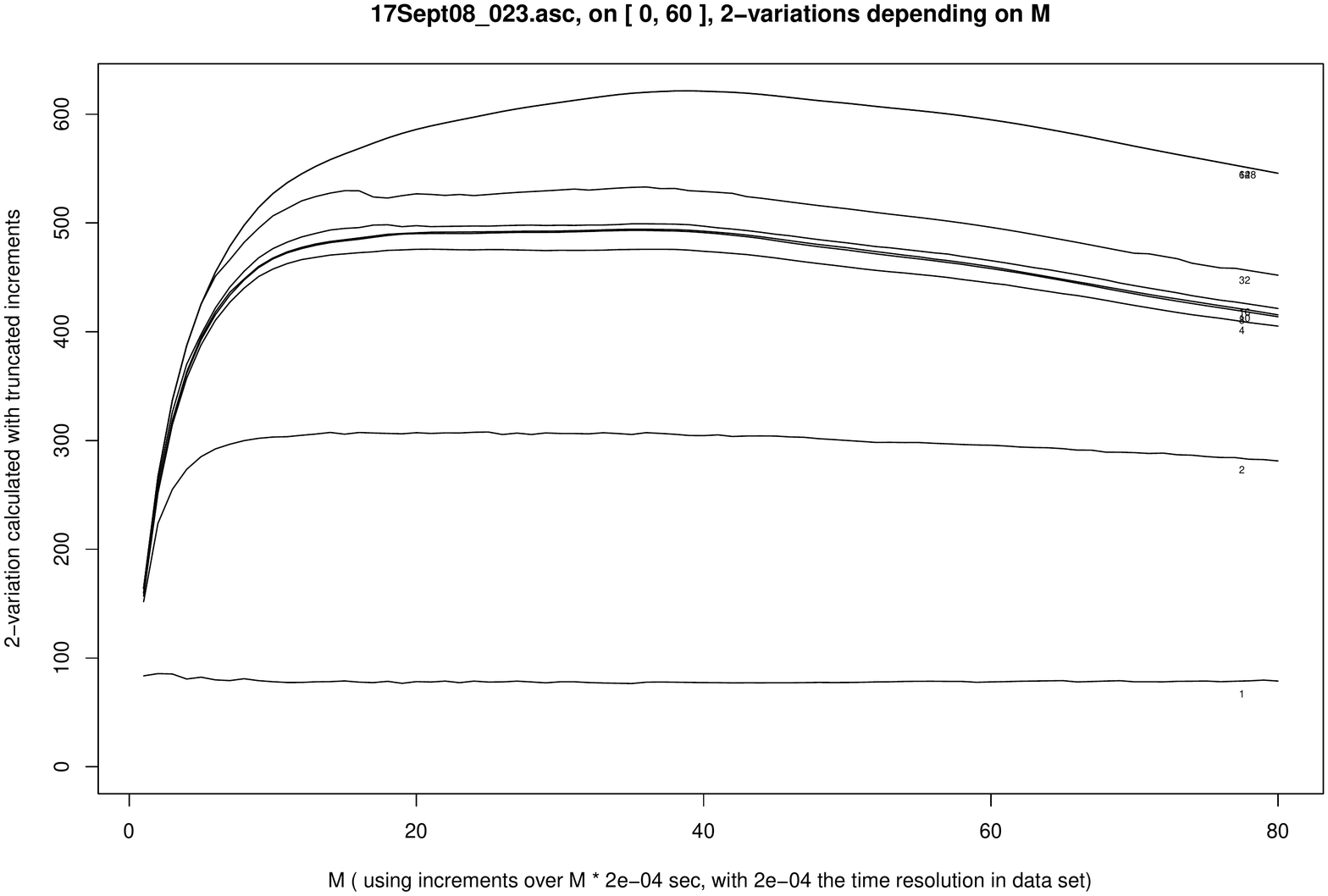}
\caption{\small Membrane potential '17Sept08{\_}023' (5 mM of K): plotting 2-variations $M\to V_\Gamma(2,\Delta,M)$ as defined in (\ref{mypvariationsegment}), with truncation constant $\Gamma$, in increasing order for $\Gamma \in\{1,2,4,8,10,16,32, 64, 128\}$; no further changes above  $\Gamma=64$.}
\label{fig10}
\includegraphics[width=12.5cm]{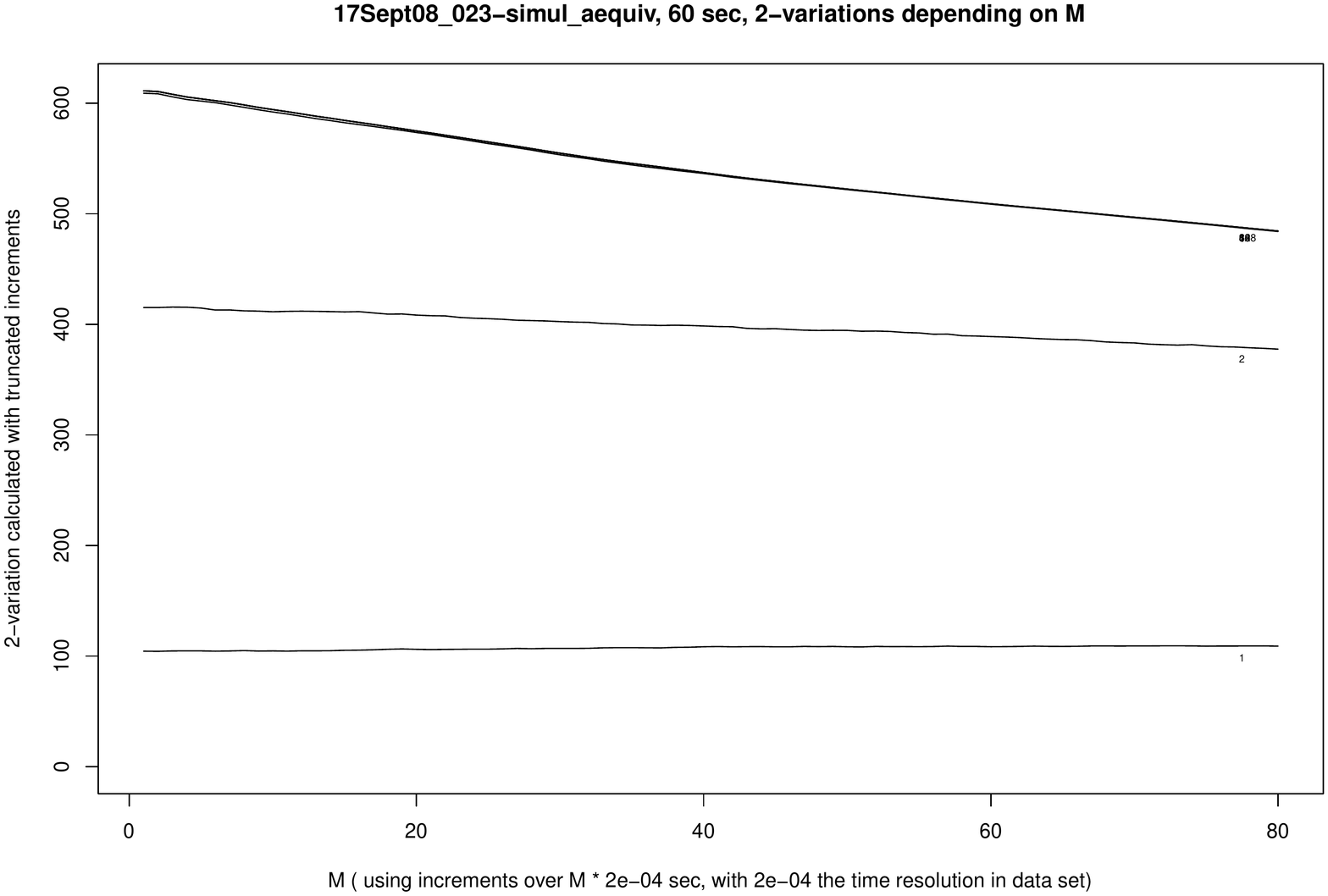} 
\caption{\small Simulated CIR type diffusion equivalent (\cite{Jahn-09}, section 5.3) for 17Sept08{\_}023 : 2-variations calculated in analogy to figure \ref{fig10}.}
\label{fig11}
\end{center}
\end{figure}

If we admit heuristics (\ref{eq-4}), then figures \ref{fig1}+\ref{fig3} ($2$-variations and $4$-variations for level 7 of 'Zelle~3', no spikes) in comparison to figures \ref{fig2}+\ref{fig4} ($2$-variations and $4$-variations for a simulated P type diffusion equivalent, with drift and diffusion coefficient as estimated in \cite{H-07}, section 4.4) show that level 7 of 'Zelle~3' exhibits the typical features of a continuous semimartingale, up to some strong deformation of the initial part of the $2$-variation as a function of $M$, visible for $M$-values up to $\approx 10$. We interpret this deformation as noise contaminating the observation, generated by the electrode measuring the membrane potential. This is supported by the observation that in all levels $1, \ldots, 10$ of 'Zelle~3' the $2$-variations attain a value close to $\approx 300$ for $M = 1$ (see figures \ref{fig9}+\ref{fig5}+\ref{fig1}), whereas in the recording of '17Sept08{\_}023' a different type of electrode was used which does not produce the same phenomenon (figure~\ref{fig10}).

Qualitative agreement (disregarding the effect of noise for $M\le 10$) between $2$- and $4$-variations for the neuronal data and $2$- and $4$-variations for their simulated diffusion equivalents is observed in all spikeless levels $1, \ldots, 7$ of 'Zelle~3', and also in level 8 where one single spike is generated in the 60 seconds of observation. This agreement is  not always as perfect as in figures \ref{fig1}+\ref{fig2} and \ref{fig3}+\ref{fig4} (presenting level 7 of 'Zelle~3' in the figures, we did chose the level where the best fit occurred), but the qualitative  features ($2$-variations flat in $M$ up to the initial effect of noise, $4$-variations linear in $M$) agree well between data and simulated diffusion equivalent. Hence in the spikeless or rarely spiking levels of 'Zelle~3', the membrane potential (away from the isolated spike in case of level 8) {\em can be be viewed as a semimartingale}. 

Is the semimartingale continuous, or does it have jumps? The $4$-variations for levels 1--8 look at first glance very much like being linear in $M$ for small $M$-values. In some of these levels however, a closer look to small $M$ values might suggest presence of jumps. As an example, comparing figures \ref{fig3}+\ref{fig4} (level 7) for $M=3,2,1$, there is a difference in the behaviour for small $M$, and in figure \ref{fig3} we may see convergence to some strictly positive 'limit' as $M$ gets small. Thus there might be jumps in the membrane potential data 'Zelle~3' level 7. The same effect is visible e.g.\  in levels 1 (3 mM of K), 3 (5 mM of K), 4 (6 mM of K).  However, since in all pictures of $2$-variations we saw additional variation for small values of $M$, with the interpretation of noise of the measuring electrode, this noise might similiarly affect the $4$-variations for small values of $M$. Hence, with the methods of section~1, we cannot decide whether or not the non- or rarely spiking levels of 'Zelle~3' should be viewed as continuous semimartingales or as semimartingales with jumps. We will be able to answer this question with a different method in section~2. 


\begin{figure}[p]
\begin{center}
\includegraphics[width=12.5cm]{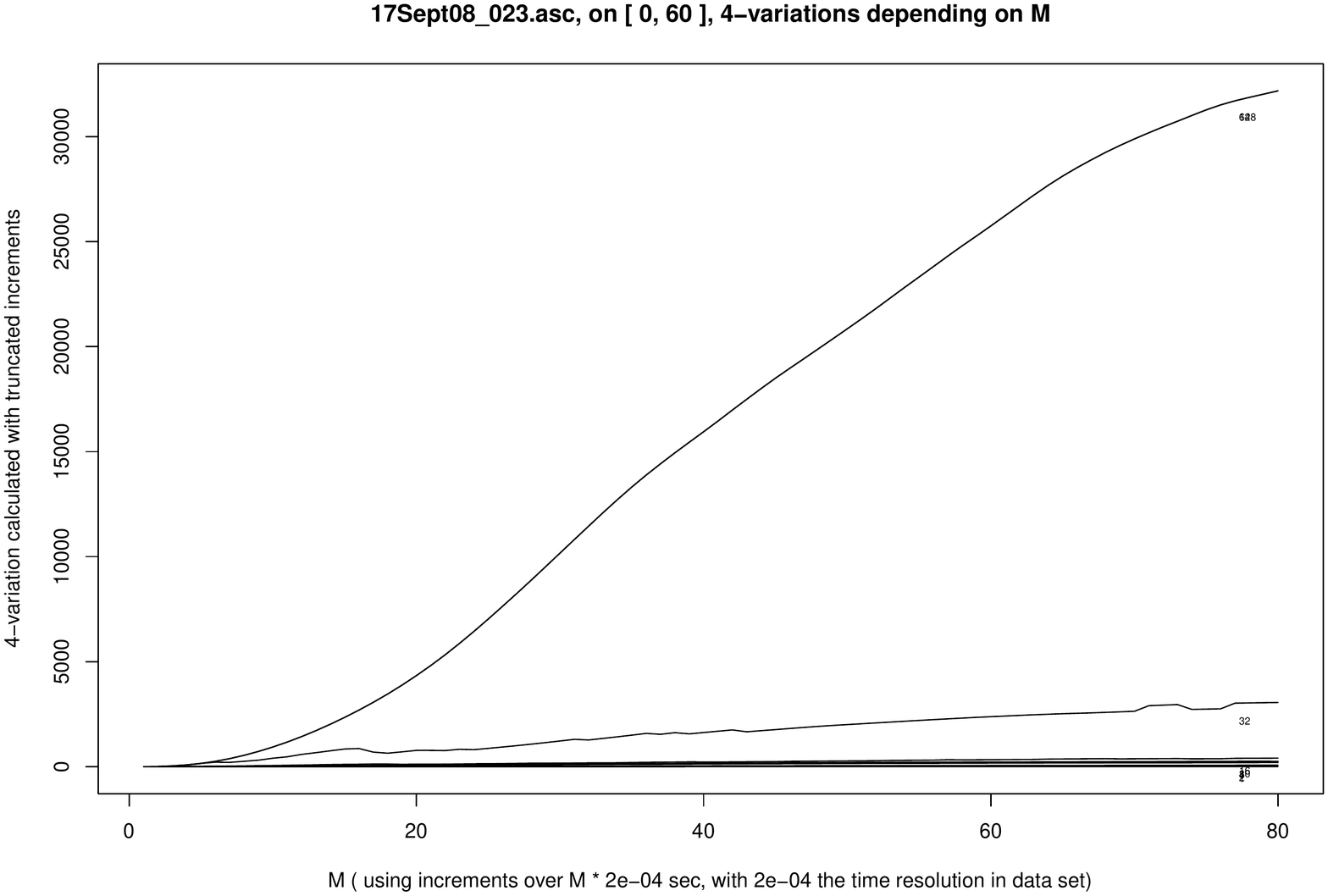}
\caption{\small Membrane potential '19Sept08{\_}023' (5 mM of K): plotting 4-variations $M\to V_\Gamma(4,\Delta,M)$ as defined in (\ref{mypvariationsegment}) with truncation constant $\Gamma$, in increasing order for $\Gamma \in\{1,2,4,8,10,16,32, 64, 128\}$; no changes for $\Gamma\ge 64$. }
\label{fig12}
\includegraphics[width=12.5cm]{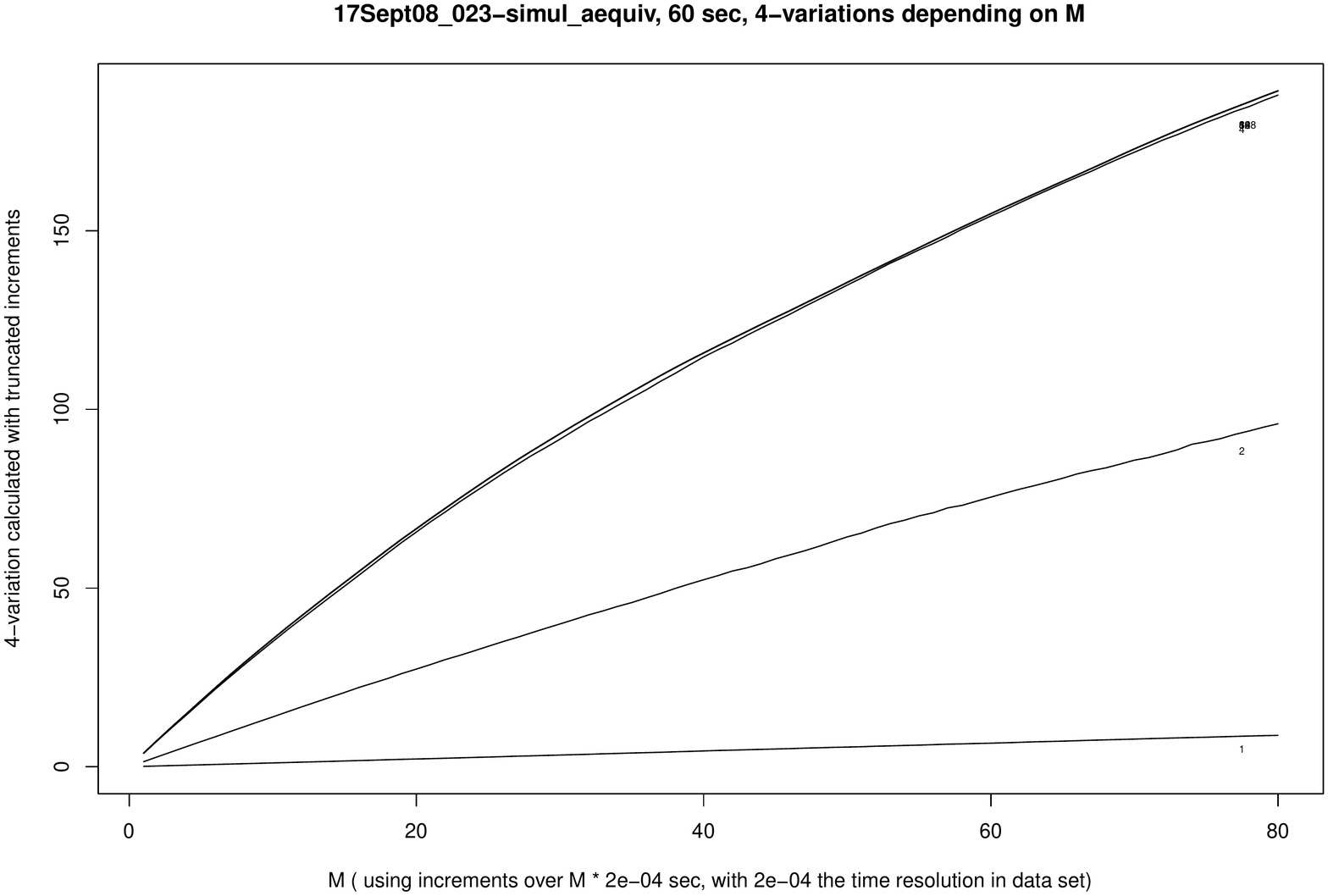} 
\caption{\small Simulated CIR type diffusion equivalent (\cite{Jahn-09}, section 5.3) for '17Sept08{\_}023': 4-variations calculated in analogy to figure \ref{fig12}.}
\label{fig13}
\end{center}
\end{figure}

\subsection*{1.3. Application to the spiking levels of 'Zelle~3', and to the neuron '17Sept08{\_}023'} 
In the spiking levels 9 and 10 of 'Zelle~3' (18 spikes in level 9, 8 spikes in level 10, over a total of 60 seconds of observation time) and in the rapidly spiking neuron '17Sept08{\_}023' ($\approx 50$ spikes over the first 60 seconds), pictures very different from those discussed above arise (see figures \ref{fig10}+\ref{fig12} in comparison to \ref{fig11}+\ref{fig13} for '17Sept08{\_}023', see figures \ref{fig5}+\ref{fig7} in comparison to \ref{fig6}+\ref{fig8} for 'Zelle~3' level 9, see also figure \ref{fig9} for 'Zelle~3' level 10). There is no longer a qualitative coincidence between the shape of the $p$-variations as a function of $M$ in the data, and the shape of the $p$-variations as a function of $M$ in the simulated diffusion equivalent. 

\begin{figure}[t]
\begin{center}
\includegraphics[width=12.5cm]{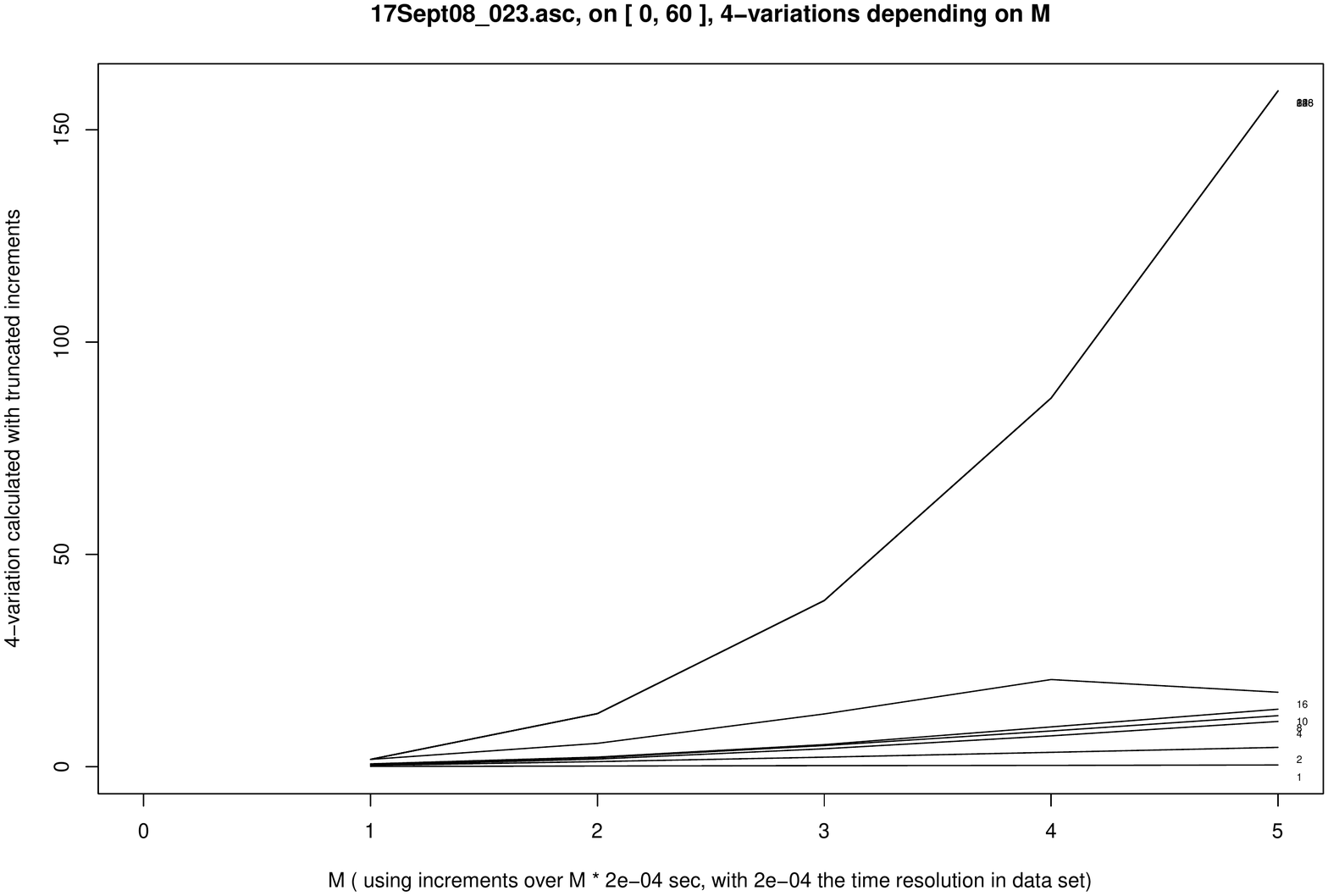}
\caption{\small Zooming into figure \ref{fig12} for small values of $M$: for $\Gamma\ge 64$, we obtain the values $39.13$ for $M=3$, $12.50$ for $M=2$, $1.73$ for $M=1$. }
\label{fig14}
\end{center}
\end{figure}

For the $4$-variations of '17Sept08{\_}023' (where thanks to $\Delta = 2\cdot 10^{-4}$ [sec] we get nearer to $0$ than in the data 'Zelle~3') we calculate the values $39.13$ at $M=3$, $\,12.50$ at $M=2$, $\,1.73$ at $M=1$ for sufficiently large truncation factor $\Gamma$, see the detail in figure \ref{fig14}. Under $p=4$, these values rule out at the same time 
\beao
&&\bullet\quad\mbox{the possibility of a strictly positive 'limit' in (\ref{eq-5neu}) when $M$ gets small} \\
&&\bullet\quad\mbox{the possibility of a linear dependence on $M$ in  (\ref{eq-7neu}) when $M$ gets small}
\eeao
and thus --according to the dichotomy in (\cite{AJ-09a}, theorem 1)-- rule out the possibility that the data '17Sept08{\_}023' represent a discretely observed semimartingale $\xi=(\xi_t)_{t\ge 0}$. With respect to this result, problems such as time inhomogeneity, obviously present in at least a part of our data, or presence of jumps as raised in \cite{H-07} become irrelevant. 
Curves of similiar shape are obtained for the $4$-variations in 'Zelle~3' level 9 (figure \ref{fig7}) and level 10. 
The situation for small values of $M$ is less clear in the data 'Zelle~3' than in '17Sept08{\_}023' (the $4$-variations in 'Zelle~3' level 9 take the values $8.17$ at $M=3$, $\,5.42$ at $M=2$, $\,3.20$ at $M=1$, cf.\ figure \ref{fig7}, and may be affected by the noise which appears clearly in the $2$-variations for small $M$). In all three data sets, the shape of the $2$-variations as function of the step size (figures \ref{fig5}, \ref{fig9}, \ref{fig10}), with a remarkable maximum at $M\approx 80$ for 'Zelle~3', at $M\approx 40$ for  '17Sept08{\_}023', is incompatible with a discretely observed semimartingale, cf.\ (\ref{eq-6}). 
Thus we conclude that the membrane potential in the spiking neuron --between successive spikes, staying sufficiently away from the spikes-- is {\em not adequately modelled by a semimartingale}. 

We remark that despite this fact, nonparametric estimates for drift and diffusion coefficient --within a semimartingale setting-- may produce seemingly satisfactory results. As an example, for 'Zelle~3' level 10, one obtains a convincing fit between occupation time calculated from the data and the invariant Gamma type law of the CIR type diffusion equivalent for 'Zelle~3' level 10 (with estimated drift and diffusion coefficient according to \cite{H-07}, section 3.2). But also here, the hypothesis of a diffusion process became questionable when in case of 'Zelle~3' level 10 the estimates used in \cite{H-07} were observed to depend much more on the chosen multiple $M$ of the step size $\Delta$ --entering the definition of the kernel estimator in \cite{H-07} -- than was claimed in \cite{H-07}. This observation represents a surprising contrast to what has been checked for the non-spiking levels 3, 6, 7 of 'Zelle~3' in (\cite{H-07}, figure 10).

\subsection*{1.4. An additional remark} 
In all levels of 'Zelle~3', the $2$-variations $M\to V_\Gamma(2,\Delta,M)$ present periodic deformations,  for fixed value of the truncation factor $\Gamma$; to less extent, this is visible also in the $4$-variations $M\to V_\Gamma(4,\Delta,M)$ (see figures \ref{fig1}+\ref{fig3}, \ref{fig5}+\ref{fig7}, \ref{fig9}).  
Independently of the level and of the value of $\Gamma$, the deformations are most visible near $M\approx 32$ and $M\approx 64$ (in figure \ref{fig9} continued by $\approx 96, 128, \ldots$), and go attenuating as $M$ gets larger. This might indicate that the neuronal network in the slice to which 'Zelle~3' belongs possesses loops or circuits, and thus produces feedback at a fixed periodicity.

\section*{2. Fixing the size of the increments and varying the power $p$} 
In section 1, we have considered $p$-variations for fixed $p$ as a function of the step size. We continue this discussion, but now with $p$-variations considered as a function of $p$ for fixed step size. 
\begin{figure}[t]
\begin{center}
\includegraphics[width=12.5cm]{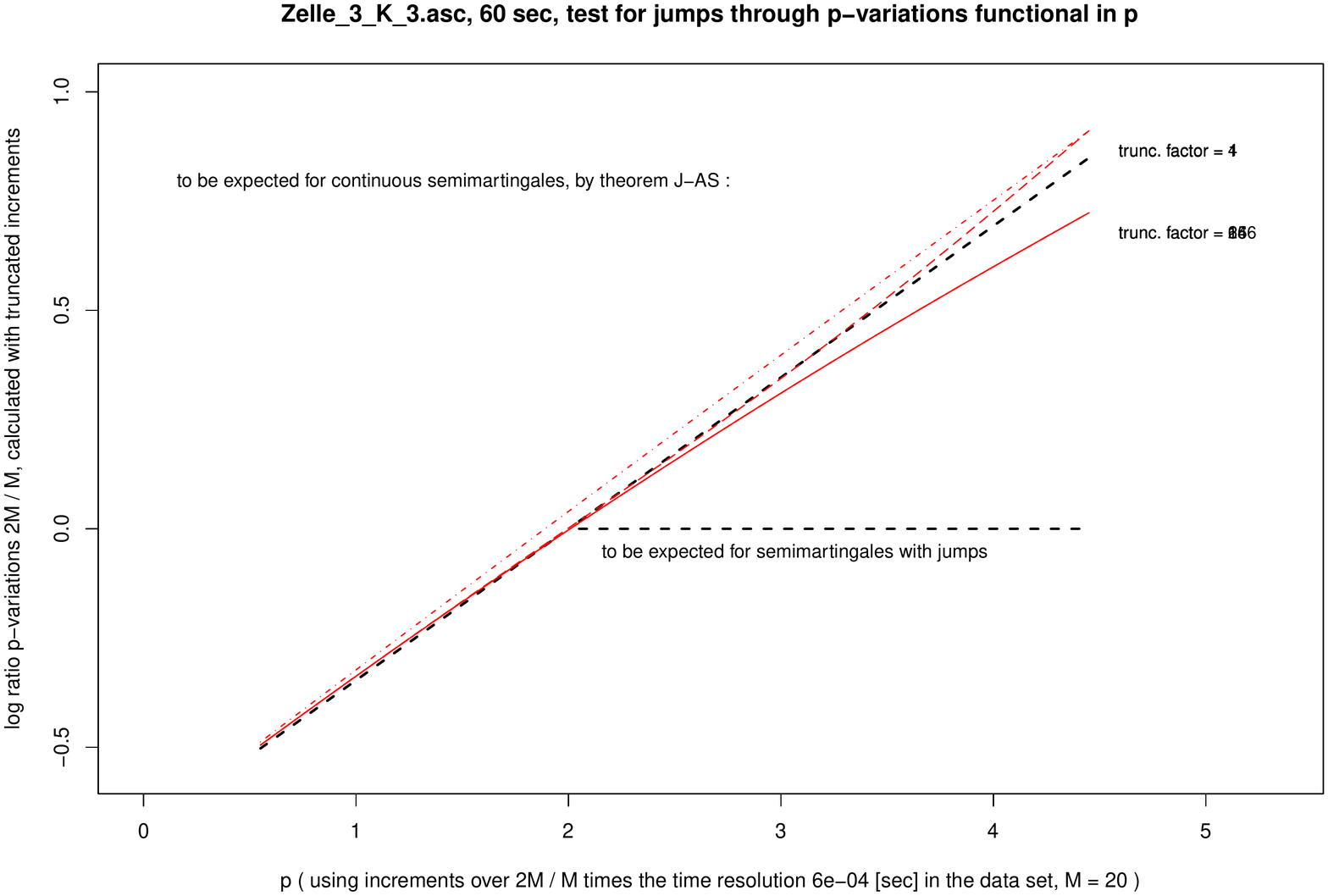} 
\caption{\small Membrane potential 'Zelle~3' level 1 (3 mM of K, no spikes). Logarithm of ratios $\frac{V_\Gamma(p,\Delta,2M)}{V_\Gamma(p,\Delta,M)}$ as a function of $p$, for $\Gamma\in\{1,4,16,64,256\}$. Increasing values of the truncation factor correspond to more 'solid' red curves; no changes occur above $\Gamma=16$. It is seen that the data 'Zelle~3' level 1 are well compatible with a continuous semimartingale observed at discrete times $i\Delta$, $0\le i\le 100001$.}
\label{fig15}
\end{center}
\end{figure}
\begin{figure}[t]
\begin{center}
\includegraphics[width=12.5cm]{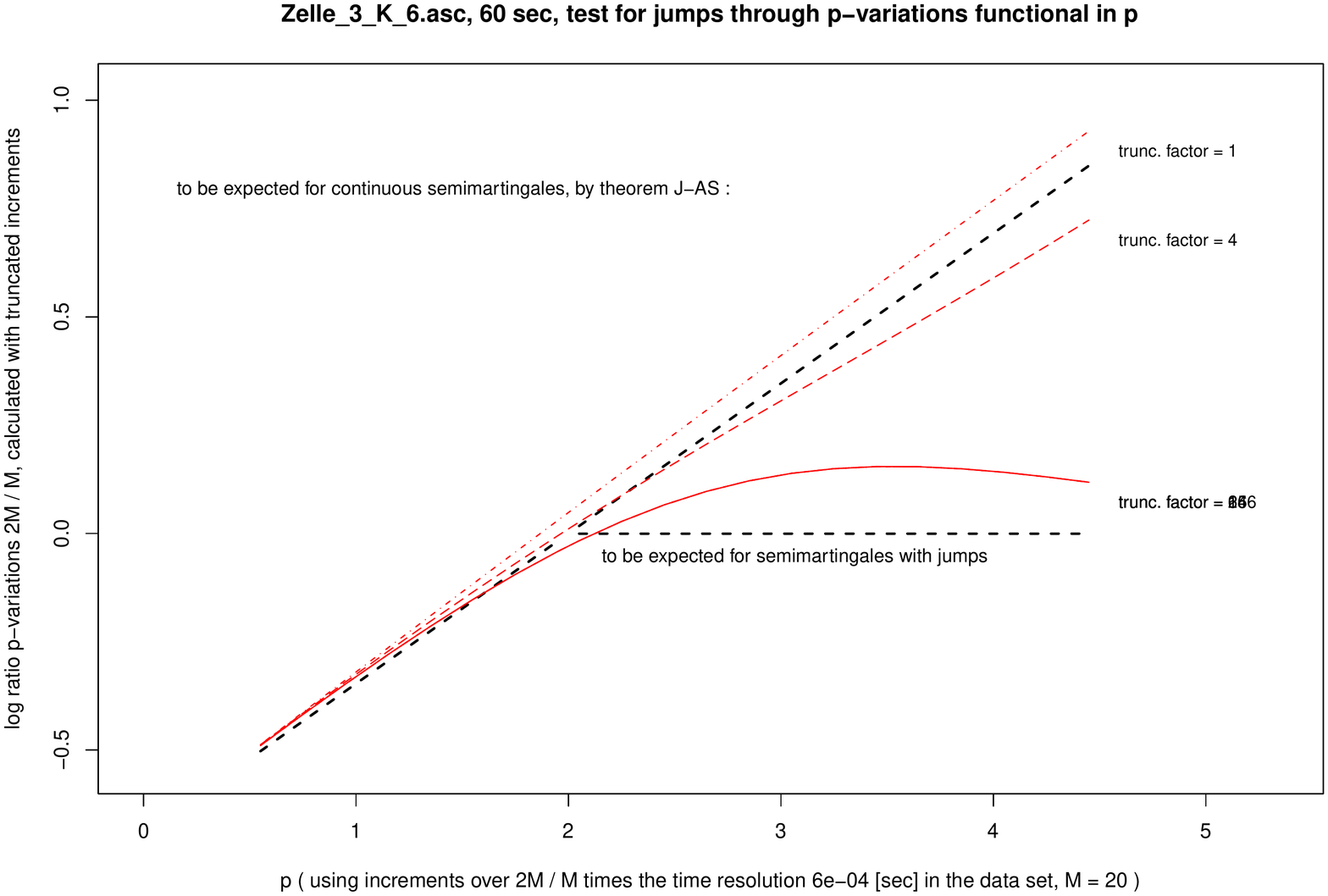}
\caption{\small Membrane potential 'Zelle~3' level 4 (6 mM of K, no spikes): logarithm of ratios $\frac{V_\Gamma(p,\Delta,2M)}{V_\Gamma(p,\Delta,M)}$ as a function of $p$, for $\Gamma\in\{1,4,16,64,256\}$. Increasing values of the truncation factor correspond to more 'solid' red curves, no changes above $\Gamma=16$. It is seen that the data 'Zelle~3' level 1 are well compatible with a semimartingale which has  jumps, observed at discrete times $i\Delta$, $0\le i\le 100001$.}
\label{fig16}
\end{center}
\end{figure}
Assuming that the spikeless segments of the membrane potential do correspond to a discretely observed semimartingale with nonvanishing continuous local martingale part (an unproblematic additional assumption given the shape of the power variations for small truncation factor $\Gamma$ appearing in figures \ref{fig1}, \ref{fig3}, \ref{fig5}, \ref{fig7}, \ref{fig10}, \ref{fig12}, \ref{fig14} above) and that both $M'\in \{ M , 2 M \}$ lead to sufficiently small values of $M'\Delta$, we may read (\cite{AJ-09a}, (11)+(7)+(10)) on a segment $[t_0,t_1]$ as follows: 
\beam\label{eq-8}
\mbox{for $\xi$ with jumps} &:& 
\wh B_{t_0,t_1}(p,\Delta,2M) \;\approx\; 
\left\{\begin{array}{ll}
\wh B_{t_0,t_1}(p,\Delta,M)      &  \mbox{for $2\le p<\infty$} \\
2^{\frac{p}{2} -1}\, \wh B_{t_0,t_1}(p,\Delta,M)    &  \mbox{for $0<p<2$}  
\end{array}\right.\\
\label{eq-9}
\mbox{for $\xi$ continuous} &:&   \wh B_{t_0,t_1}(p,\Delta,2M) \quad\approx\quad 2^{\frac{p}{2} -1}\, \wh B_{t_0,t_1}(p,\Delta,M) \;\;\mbox{for $0<p<\infty$}  
\eeam
Again we accept the heuristics of section 1.1, in particular the approximation (\ref{eq-4}). We extend the heuristics by assuming that if spikeless segments of the membrane potential do correspond to a semimartingale which has jumps, then jumps will occur on every segment $[t_0,t_1]$ under consideration (this is unproblematic e.g.\ if the L\'evy measure of the jump part of $\xi$ has infinite total mass independently of time). Then we can rephrase the test for jumps in Ait-Sahalia and Jacod (\cite{AJ-09a}, theorem 1) for fixed $M$ and varying $p$ as follows: 
\beam\label{eq-10}
\mbox{for $\xi$ with jumps}&:&  p \lra \log\frac{V(p,\Delta,2M)}{V(p,\Delta,M)} \;\;\mbox{is approximately} 
\left\{\begin{array}{l} 
\mbox{constant  $\equiv 0\;$ on $[2,\infty)$} \\
\mbox{linear in $p\;$ on $(0,2)$}  
\end{array}\right. \\ 
\label{eq-11}
\mbox{for $\xi$ continuous}&:&    p \lra \log\frac{V(p,\Delta,2M)}{V(p,\Delta,M)} \;\;\mbox{is approximately linear in $p\;$ on $(0,\infty)$} \;. 
\eeam
The slope of the linear parts in (\ref{eq-10})+(\ref{eq-11}) is deterministic, by (\ref{eq-8})+(\ref{eq-9}). Thus, from the very beginning, we know the shape which we expect to see when the membrane potential data between successive spikes do correspond to a discretely observed semimartingale $(\xi_t)_{t\ge 0}$:  in this case, the empirical object 
$$ 
p \;\lra\; \log\frac{V(p,\Delta,2M)}{V(p,\Delta,M)}  
$$
in (\ref{eq-10})+(\ref{eq-11}) should be close to 
\beam\label{eq-10neu}
p &\lra& \min \left\{ (\frac{p}{2}-1)\log 2 \,,\, 0 \right\}  \quad  \mbox{on $\;0<p<\infty$} \quad \mbox{if $\xi$ has jumps} \;,  
\\ \label{eq-11neu}
p &\lra& (\frac{p}{2}-1)\log 2  \quad  \mbox{on $\;0<p<\infty$} \quad \mbox{if $\xi$ is continuous} \;.  
\eeam
Comparing the empirical object in (\ref{eq-10})+(\ref{eq-11}) to the truncated line (\ref{eq-10neu}) expected for a semimartingale having jumps, or to the straight line (\ref{eq-11neu}) expected for a continuous semimartingale, we can decide whether jumps are present. Beyond this, since the above (\ref{eq-10})+(\ref{eq-11}) represents a dichotomy on a very general class of semimartingales, see \cite{AJ-09a}, we can decide whether or not our membrane potential data (away from the spikes) do correspond to a semimartingale.

When visualizing the empirical object in (\ref{eq-10})+(\ref{eq-11}), we will continue to make use of the truncation factor $\Gamma$ as in (\ref{mypvariation})+(\ref{mypvariationsegment}), and calculate from our data log-ratios  
\beqq\label{mypvariationratiologarithmic}
p \;\lra\; \log\frac{V_\Gamma(p,\Delta,2M)}{V_\Gamma(p,\Delta,M)}  \;. 
\eeqq
for varying values of $\Gamma$ which are representative for asymptotics $\Gamma\to\infty$. 

\begin{figure}[t]
\begin{center}
\includegraphics[width=12.5cm]{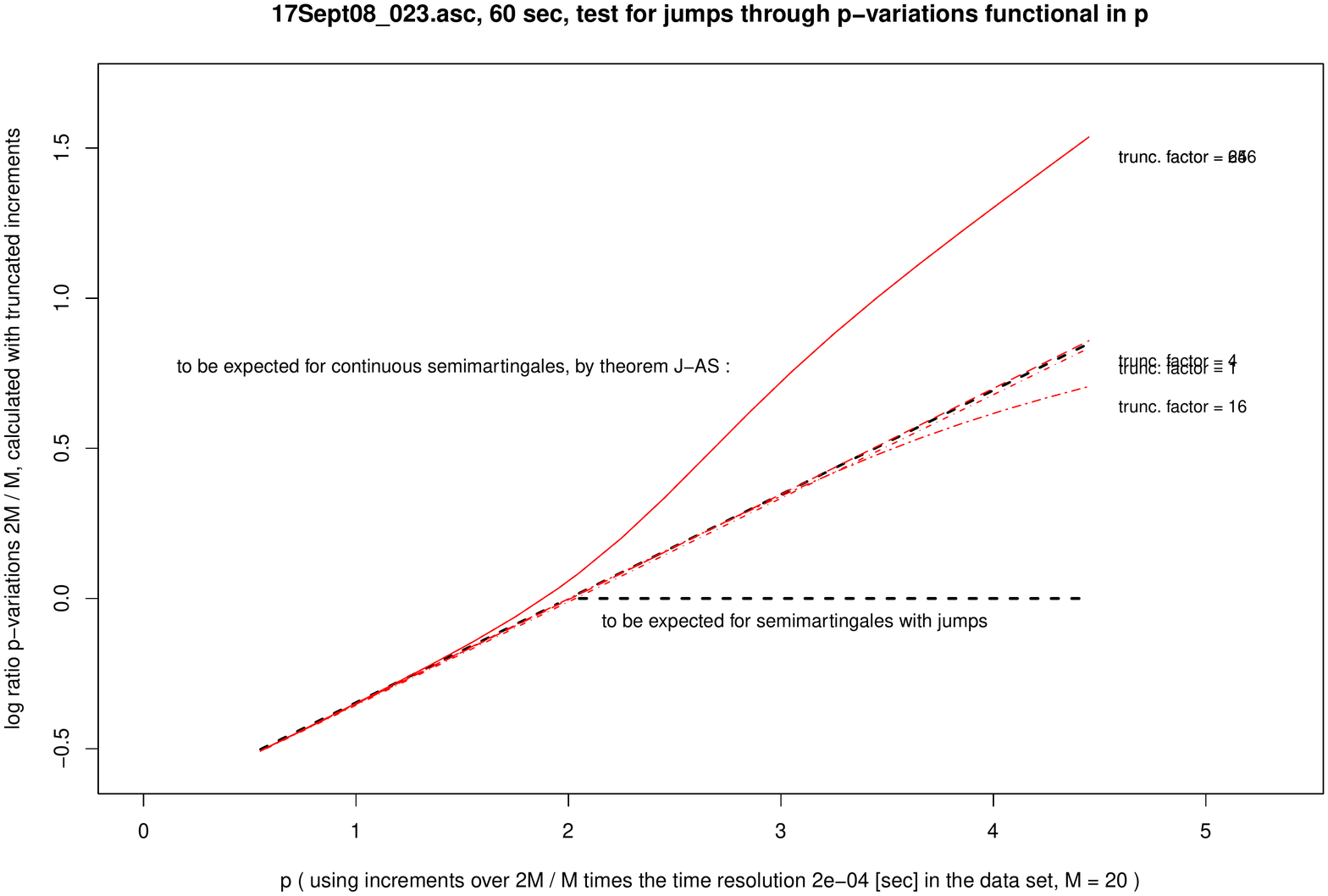}
\caption{\small Frequently spiking neuron '17Sept08{\_}023' (5 mM of K, $\approx 50$ spikes over the first 60 seconds): logarithm of ratios $\frac{V_\Gamma(p,\Delta,2M)}{V_\Gamma(p,\Delta,M)}$ plotted as a function of $p$, for $\Gamma\in\{1,4,16,64,256\}$; the curves stabilize for $\Gamma\ge 64$.  Increasing values of the truncation factor correspond to more 'solid' red curves. The dotted black lines show what is expected for a semimartingale, continuous or not, by \cite{AJ-09a}. For $\Gamma$ tending to $\infty$, the logarithmic ratios calculated from the data '17Sept08{\_}023' turn out to be quite far away from a semimartingale hypothesis.}
\label{fig17}
\end{center}
\end{figure}

\begin{figure}[t]
\begin{center}
\includegraphics[width=12.5cm]{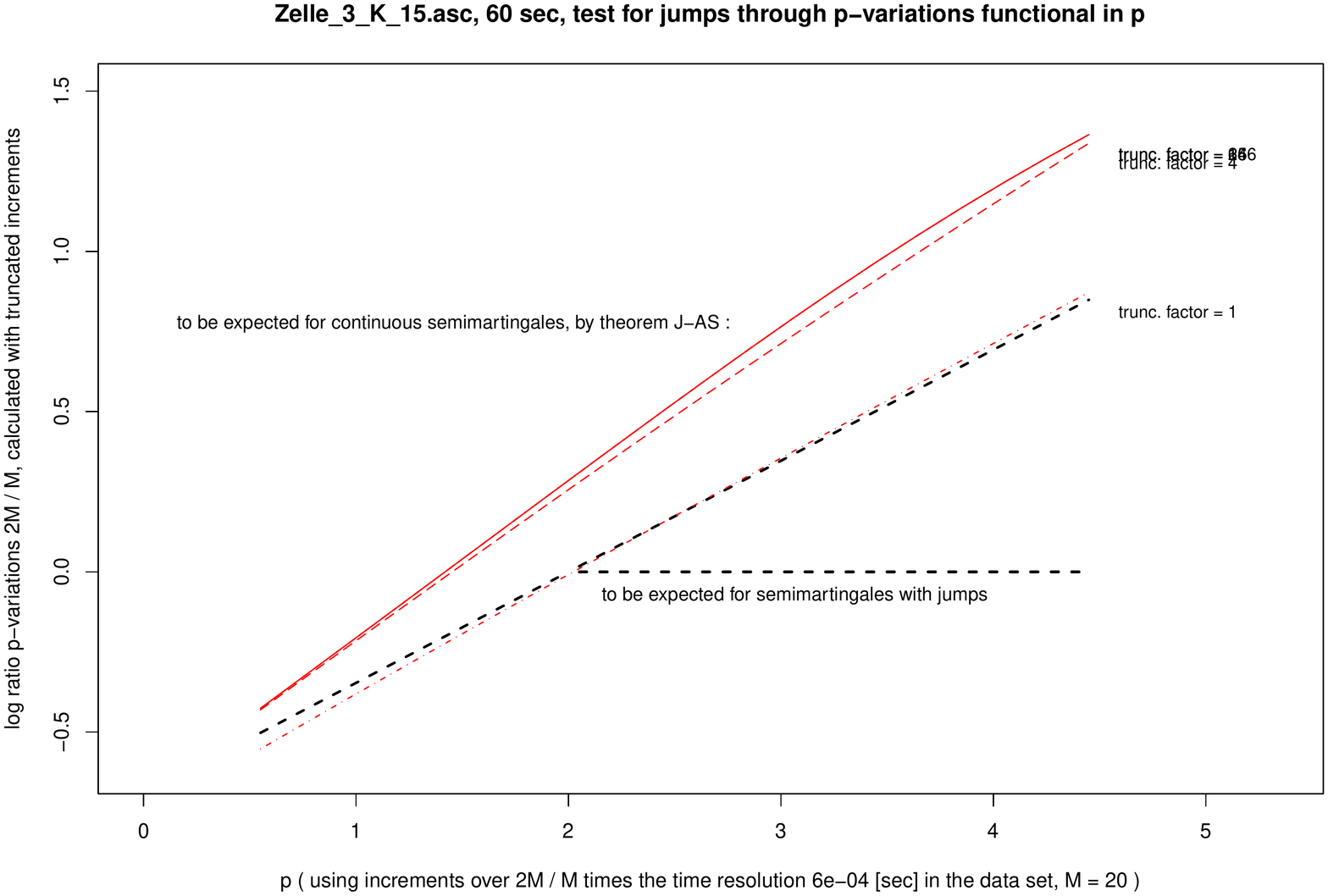} 
\caption{\small 'Zelle~3' level 10 (15 mM of K, 8 spikes over 60 seconds): logarithmic ratios $\frac{V_\Gamma(p,\Delta,2M)}{V_\Gamma(p,\Delta,M)}$ plotted as a function of $p$, for $\Gamma\in\{1,4,16,64,256\}$; the curves stabilize for $\Gamma\ge 64$.  Increasing values of the truncation factor correspond to more 'solid' red curves. The dotted black lines show what is expected for a semimartingale, continuous or not, by \cite{AJ-09a}. For $\Gamma$ tending to $\infty$, the logarithmic ratios calculated from the data '17Sept08{\_}023' do not fit well with a semimartingale hypothesis.}
\label{fig18}
\end{center}
\end{figure}

\subsection*{Application to the data 'Zelle~3' and to '17Sept08{\_}023'}
Figure \ref{fig15} shows the lowest level of 'Zelle~3' (3 mM of K, spikeless): the picture corresponds very well to what we expect for a semimartingale which is continuous, as explained in (\ref{eq-10neu})+(\ref{eq-11neu})+(\ref{mypvariationratiologarithmic}) above. Figure \ref{fig16} shows the level 4 of 'Zelle~3' (6 mM of K, spikeless): the shape of the curve corresponds very well to what we expect for a semimartingale which has jumps. Simulated diffusion equivalents (resp.: simulating a jump diffusion as in figure \ref{figeffectsofjumps}, in relation to 'Zelle~3' level 5) produce pictures similiar to figure \ref{fig15} (resp.: to figure \ref{fig16}). Moreover, {\em all non-spiking levels} 1--7 of 'Zelle~3', and even level 8 with one isolated spike over 60 seconds of observation,  lead to curves {\em corresponding convincingly to a semimartingale hypothesis}. Among these, exactly two --the levels 4 and 5 (level 5, not shown, looks much like figure \ref{fig16})-- indicate the presence of jumps. Thus the method used in the present section, in contrast to the method used in section~1, is able to answer the problem of jumps raised in section~1. 

It turned out in section~1 that the membrane potential in the spiking levels 9 and 10 of 'Zelle~3' and in the frequently spiking neuron '17Sept08{\_}023' was not adequately modelled by a semimartingale. The method of the present section reinforces this, see figures \ref{fig17}+\ref{fig18}. The picture for level 9 of 'Zelle~3' (not shown) is similiar to what we show for level 10 in figure \ref{fig18}.   The three curves do not correspond to what we expect for a semimartingale --continuous or not-- by \cite{AJ-09a}. Note that in levels 9 or 10 of 'Zelle~3', the time intervals between successive spikes are still relatively large and contain enough observations such that semimartingale methods --if the observed process were a semimartingale-- should work successfully. 

On the basis of (\ref{eq-10neu})+(\ref{eq-11neu}), we conclude as in section~1 that in the spiking levels of 'Zelle~3' and in the frequently spiking neuron '17Sept08{\_}023', the membrane potential between successive spikes {\em is not a discretely observed semimartingale}. Note that our data 'Zelle~3' have been collected {\em in the same neuron in the same cortical slice under different level of stimulation} where stimulation by potassium activates the networking properties of all neurons in the slice. The membrane potential of the observed neuron --sufficiently away from the spike times whenever there are spikes-- behaves as a semimartingale as long as there are no spikes or at most extremely isolated ones, and loses this property once spikes occur frequently enough. This adresses in particular a serious question to some widely used neuronal models where interspike intervals are identified with level crossing times of semimartingales.

\begin{figure}[p]
\begin{center}
\includegraphics[width=17.5cm]{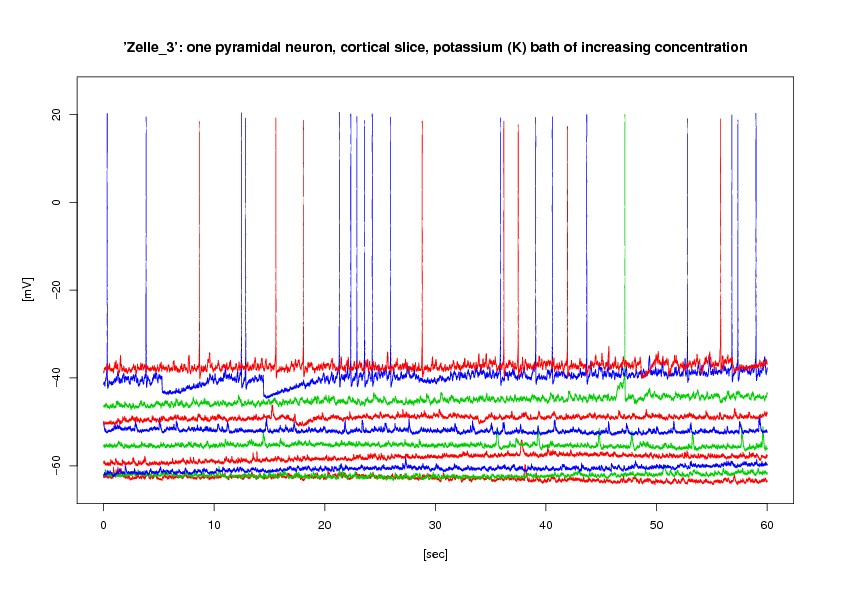}
\caption{\small 'Zelle~3', all levels 1--10: membrane potential in the same pyramidal neuron under different experimental conditions. The neuron belongs to a cortical slice observed in vitro. The networking properties of all neurons in the slice are stimulated by a  potassium bath (3, 4, 5, 6, 7, 8, 9, 10, 12, 15 mM of K). Spikes occur in levels 9 (18 spikes within 60 seconds of observation time) and level 10 (8 spikes), one isolated spike being observed in level 8. The time resolution is $\Delta = 6\cdot 10^{-3}$ [sec].  Data from H.\ Luhmann and W.\ Kilb, Institute of Physiology, University of Mainz. }
\label{fig19}
\end{center}
\end{figure}

\begin{figure}[p]
\begin{center}
\includegraphics[width=17.5cm]{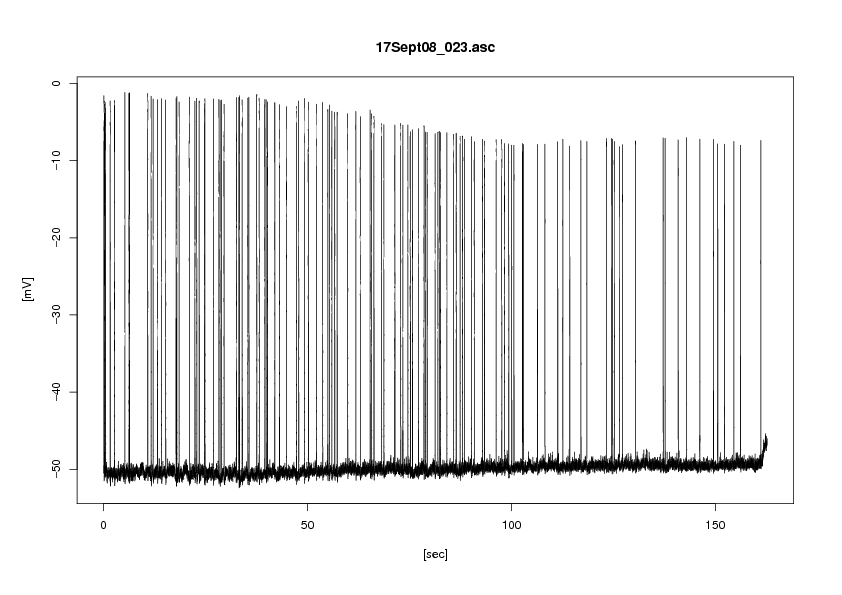}
\caption{\small Membrane potential in the frequently spiking neuron '17Sept08{\_}023' (5 mM of K); we will use only the part of the data which corresponds to the first 60 seconds of observation. On this time interval, approximately 50 spikes occur. The time resolution is $\Delta = 2\cdot 10^{-3}$ [sec]. Data from H.\ Luhmann and W.\ Kilb, Institute of Physiology, University of Mainz.}
\label{fig20}
\end{center}
\end{figure}

\vskip1.5cm 

\vskip1.5cm

~ \hfill{\bf 21.08.2010}\\

Reinhard H\"opfner, 
Institut f\"ur Mathematik, Universit\"at Mainz, 
Staudingerweg 9, D--55099 Mainz, Germany\\
{\tt hoepfner@mathematik.uni-mainz.de},  
{\tt www.mathematik.uni-mainz.de/$\sim$hoepfner}

\end{document}